\documentclass{amsart}
\usepackage{amssymb}
\usepackage{amsfonts}
\usepackage{amstext}
\usepackage{algorithmic}
\usepackage{algorithm}
\usepackage{graphicx}
\usepackage{epstopdf}
\usepackage[all]{xy}

\numberwithin{equation}{section}
\newtheorem{theorem}{Theorem}[section]
\newtheorem{lemma}[theorem]{Lemma}
\newtheorem{proposition}[theorem]{Proposition}
\newtheorem{corollary}[theorem]{Corollary}
\newtheorem{conjecture}[theorem]{Conjecture}
\newtheorem{question}[theorem]{Question}
\newtheorem{fact}[theorem]{Fact}
\theoremstyle{definition}
\newtheorem{definition}[theorem]{Definition}
\newtheorem{example}[theorem]{Example}

\theoremstyle{remark}
\newtheorem{remark}[theorem]{\bf{Remark}}

\newtheorem{notation}[theorem]{\bf{Notation}}

\newcommand{\Sym}{{\rm{Sym}}}

\newcommand{\Acal}{{\mathcal{A}}}
\newcommand{\C}{{\mathbb{C}}}
\newcommand{\Z}{{\mathbb{Z}}}
\newcommand{\N}{{\mathbb{N}}}

\newcommand{\Hcal}{{\mathcal{H}}}

\newcommand{\Fcal}{{\mathcal{F}}}
\newcommand{\Gcal}{{\mathcal{G}}}

\newcommand{\Lcal}{{\mathcal{L}}}
\newcommand{\Mcal}{{\mathcal{M}}}

\newcommand{\Pcal}{{\mathcal{P}}}

\newcommand{\Rcal}{{\mathcal{R}}}

\newcommand{\tens}{\otimes}
\newcommand{\id}{{\rm id}}
\newcommand{\bo}{{}^{(1)}}
\newcommand{\bt}{{}^{(2)}}
\renewcommand{\o}{{}_{(1)}}
\renewcommand{\t}{{}_{(2)}}

\newcommand{\extd}{{\rm d}}
\newcommand{\del}{{\partial}}

\begin{document}

\title[Multipermutation solutions of level two ]{Quantum spaces associated to multipermutation solutions of level two}
\thanks{The first author is partially supported by Isaac Newton Institute, UK, the ICTP, Trieste, and 
Grant MI 1503/2005 of the Bulgarian National Science Fund of the
Ministry of Education and Science.}

\author{Tatiana Gateva-Ivanova and Shahn Majid}
\address{TGI: Institute of Mathematics and Informatics\\
Bulgarian Academy of Sciences\\
Sofia 1113, Bulgaria\\
S.M:  Queen Mary, University of London\\
School of Mathematics, Mile End Rd, London E1 4NS, UK}

\email{tatianagateva@yahoo.com, tatyana@aubg.bg,
s.majid@qmul.ac.uk}

\date{June 2008/Ver 2}

\begin{abstract}
We study finite   set-theoretic
solutions $(X,r)$ of the Yang-Baxter equation of square-free multipermutation type. We show that each
such solution over $\C$ with multipermutation level two can be put in diagonal form with the associated Yang-Baxter algebra $\Acal(\C,X,r)$ having a $q$-commutation form of relations determined by  complex phase factors. These complex factors are roots of unity and all roots of a prescribed form appear as determined by the representation theory of the finite abelian group $\Gcal$ of left actions on $X$.  We study the structure of $\Acal(\C,X,r)$ and show that they have a $\bullet$-product form `quantizing'  the commutative algebra of polynomials in $|X|$ variables. We obtain the $\bullet$-product  both as a Drinfeld cotwist for a certain canonical 2-cocycle and as a braided-opposite product for a certain crossed $\Gcal$-module (over any field $k$). We provide first steps in the noncommutative differential geometry of $\Acal(k,X,r)$ arising from these results. As a byproduct of our work  we find that every such level 2 solution $(X,r)$ factorises as $r=f\circ\tau\circ f^{-1}$ where $\tau$ is the flip map and $(X,f)$ is another solution coming from $X$ as a crossed $\Gcal$-set.

\end{abstract}
\maketitle

\section{Introduction }

Let $V$ be a vector space of a field $k$. It is well-known that the `Yang-Baxter equation' (YBE) on an  invertible  linear map $R:V\tens V\to V\tens V$,
\[ R_{12}R_{23}R_{12}=R_{23}R_{12}R_{23}\]
(where $R_{i,j}$ denotes $R$ acting in the $i,j$ place in $V\tens V\tens V$)
provide a linear representation of the braid group on tensor powers of $V$. When $R^2=\id$ one says that the solution is involutive, and in this case
one has a representation of the symmetric group on tensor powers. We say that $R$ is of cotwist form if $R=FF_{21}^{-1}$ for some invertible linear map $F:V\tens V\to V\tens V$. Such a form of $R$ always obeys the YBE and is always involutive. Associated to each solution of the YBE are several algebraic constructions, among them the `quantum space' $k_R[V]$ which is nothing but the tensor algebra on $V$ modulo the ideal generated by the image of $\id-R$. If $\{x_i\}$ are a basis of $V$ then $k_R[V]$ is the free associative algebra generates by the $x_i$ modulo the ideal generated by  $x_ix_j-\cdot R(x_i\tens x_j)$ where $\cdot$ on the right denotes product in the free associative algebra. 

Next, a particularly nice class of solutions is provided by set-theoretic solutions where $X$ is a set, $r:X\times X\to X\times X$ obeys similar relations on $X\times X\times X$. Of course, each such solution extends linearly to $V=kX$ with matrices in this natural basis having only entries from 0,1 and many other nice properties. In this case there is a semigroup $S(X,r)$ generated by $X$ with relations $xy=\cdot r(x,y)$ (where $\cdot$ denotes product in the free semigroup) and  the semigroup algebra $k S(X,r)= k_R[V]$. There is a canonical action of $S(X,r)$ on $X$ by permutations and hence of $k_R[V]$ on $V$. The algebra $k_R[V]$ when defined in the basis $X$ is the `Yang-Baxter algebra' $\Acal(k,X,r)$ associated to $(X,r)$, but $k_R[V]$ is defined independently of any basis of $V$ and may be computed in some other more convenient one.

In this paper we look at the special properties of linear solutions obtained in this way when $r$ is a finite square-free multipermutation solution of level $2$. Precise details of what this means will be recalled in the preliminaries but, briefly, such solutions are involutive  and have the property that the permutation group $\Gcal(X,r)$, defined as the image in $\Sym(X)$  of $S(X,r)$ acting by permutations on $X$, is abelian. We show over the field $k=\C$ that then  (i) the linear extension has a cotwist form and (ii)  there exists a basis $Y$ of $V$ in which $R$ has a diagonal form in the sense that $R$ is the vector space flip map but with further coefficients in this basis. This is achieved (Theorem\ref{maintheorem}) in Section~4 using Fourier theory on the finite groups $\Gcal_i$ which are the restrictions of $\Gcal$ to its orbits in $X$. The necessary methods in the context of finite abelian groups are first developed in Section~3 before applying to each $\Gcal_i$.

The basis $Y=\{y^i_\eta\}$ which we obtain here for the span of each orbit $X_i$ is labelled by characters on $\Gcal_i$. Writing the diagonal coefficients of $R$ in this basis as $q_{i,\eta,j,\zeta}$, say,  the quantum algebra $\Acal(\C,X,r)$ takes a `q-commutation' relations form 
\[ y^i_\eta y^j_\zeta=q_{i,\eta,j,\zeta} y^j_\zeta y^i_\eta;\quad q_{i,\eta,j,\zeta}={\chi_\zeta(\sigma^i_j)\over \chi_\eta(\sigma^j_i)}\]
where $\sigma^j_i\in\Sym(X_i)$ is a certain fixed permutation for each ordered pair $i,j$ and $\chi_\eta,\chi_\zeta$ denote characters on $\Gcal_i,\Gcal_j$. Relations with $q$-factors are typical in the simplest models of noncommutative geometry (for example the `quantum plane' \cite{Man}) but now the factors are particular roots of unity with values determined by the group characters. 

Such factors are also typical in quantum mechanics (the `Weyl form' of the Heisenberg commutation relations) and in Section~6 we provide a general picture of $\Acal(k,X,r)$ as `quantisation' in the sense of a deformed $\bullet$-product on the vector space of the commutative algebra $k[V]$ (the symmetric algebra on $V$ which we view as polynomials in any fixed basis of $V$ as $|X|$ variables). We obtain this $\bullet$ product in analogy with the Moyal product in quantum mechanics, using modern Drinfeld cotwist methods. The form of $R$ ensures that this is possible but we find a natural expression of $\bullet$ in Proposition~\ref{Acotwist} as such a cotwist with respect to a canonical cocycle on $\Hcal^*\times \Hcal$ where $\Hcal=\prod_i\Gcal_i$. This is a finite group version of the canonical symplectic structure on $V^*\oplus V$ in quantum theory. We then find a second expression (now over a general field $k$) of the same $\bullet$-product in terms of the braided-opposite algebra to $k[V]$ in the braided category of crossed $\Gcal$-modules, Proposition~\ref{Abraop}.

An immediate consequence of our work is that $\Acal(k,X,r)$ is not only an algebra but a Hopf algebra in a certain symmetric monoidal category with braiding induced by $R$. Moreover, the infinitesimal translations afforded by this braided coproduct provide braided differential operators or `vector fields' which define a noncommtuative differential structure on the algebra. We show, equivalently, that there is a natural exterior algebra of differential forms $(\Omega(\Acal(k,X,r)),\extd)$ obtained by extending $\bullet$ to differential forms. In this way we provide the first stage of the noncommutative differential geometry of $\Acal(k,X,r)$ for the class of $(X,r)$ under study.

In Section~7 we end with a short epilogue to our work. Namely, the `cotwist' form of $R$ turns out to have a set-theoretic analogue $r=f\circ\tau\circ f^{-1}$ where $f:X\times X\to X\times X$ is itself a set-theoretic solution $(X,f)$ and $\tau(x,y)=(y,x)$. It is the (non-involutive) braiding associated to $X$ as a crossed $\Gcal$-set. We provide (Proposition~\ref{rcotwist}) necessary and sufficient conditions for this which include (and do not allow much more than) $(X,r)$ square-free involutive of multipermutation level $2$.

\section{Preliminaries on set-theoretic solutions}

There are many  works on
set-theoretic solutions and related structures, of which a relevant
selection for the interested reader is \cite{ESS,TM, T04, T04s,  TSh07, TSh08, GJO,JO,LYZ,Rump,
Tak, V}. 
In this section we recall basic notions and results which will be used in the paper.
 We shall use the terminology,  notation and some results from our previous works,
\cite{T04,  TSh07, TSh08}.

 \begin{definition}
Let $X$ be a nonempty set (not necessarily finite) and let $r: X\times X \longrightarrow
X\times X$ be a bijective map. We refer
to it as a \emph{quadratic set}, and denote it by  $(X,r)$. The image of $(x,y)$ under $r$ is
presented as
\begin{equation}
\label{r} r(x,y)=({}^xy,x^{y}).
\end{equation}
The formula (\ref{r}) defines a ``left action" $\Lcal: X\times X
\longrightarrow X,$ and a ``right action" $\Rcal: X\times X
\longrightarrow X,$ on $X$ as:
\begin{equation}
\label{LcalRcal} \Lcal_x(y)={}^xy, \quad \Rcal_y(x)= x^{y},
\end{equation}
for all $ x, y \in X.$ The map  $r$ is \emph{nondegenerate}, if
the maps $\Rcal_x$ and $\Lcal_x$ are bijective for each $x\in X$.
In this paper we shall consider only  the case where $r$
is nondegenerate. As a notational tool, 
we  shall sometimes identify
the sets $X\times X$ and $X^2,$ the set of all monomials of length
two in the free semigroup $\langle X\rangle.$
\end{definition}
\begin{definition} 
\label{defvariousr}
\begin{enumerate}
 \item
$r$ is \emph{square-free} if $r(x,x)=(x,x)$ for all $x\in X.$
\item \label{YBE} $r$ is \emph{a set-theoretic solution of the
Yang-Baxter equation} or, shortly \emph{a solution}
 (YBE) if  the braid relation
\[r^{12}r^{23}r^{12} = r^{23}r^{12}r^{23}\]
holds in $X\times X\times X,$ where the two bijective maps
$r^{ii+1}: X^3 \longrightarrow X^3$, $ 1 \leq i \leq 2$ are
defined as  $r^{12} = r\times \id_X$, and $r^{23}=\id_X\times r$. In
this case we shall refer to  $(X,r)$ also as a  \emph{a braided set}.
\item
A braided set $(X,r)$ with $r$ involutive is called \emph{a
symmetric set}.
\end{enumerate}
\end{definition}
To each quadratic set $(X,r)$  we associate canonical algebraic
objects  generated by $X$ and with quadratic defining relations
$\Re =\Re(r)$ defined by
\begin{equation}
\label{defrelations} xy=zt \in \Re(r),\quad
  \text{whenever}\quad r(x,y) = (z,t).
\end{equation}
\begin{definition}
\label{associatedobjects}  Let $(X,r)$ be a quadratic set.

(i) The semigroup
$ S =S(X, r) = \langle X; \Re(r) \rangle$, with a
set of generators $X$ and a set of defining relations $ \Re(r),$
is called \emph{the semigroup associated with $(X, r)$}.

(ii) The \emph{group $G=G(X, r)$ associated with} $(X, r)$ is
defined as
$G=G(X, r)={}_{gr} \langle X; \Re (r) \rangle$

(iii) For arbitrary fixed field $k$, \emph{the $k$-algebra
associated with} $(X ,r)$ is defined as
$\Acal = \Acal(k,X,r) = k\langle X ; \Re(r) \rangle.$ 
($\Acal(k,X,r)$ is isomorphic to the monoidal algebra $kS(X,r)$).

(iv) To each nondegenerate braided set $(X,r)$ we also associate a permutation group, called
\emph{the group of left action} and denoted $\Gcal= \Gcal(X,r)$, see Definition \ref{Gcaldef}.

If $(X,r)$ is a solution, then $S(X,r)$, resp. $G(X,r)$, resp.
$\Acal(k,X,r)$ is called the {\em{Yang-Baxter semigroup}}, resp.
the {\em{Yang-Baxter group}}, resp. the {\em{Yang-Baxter
algebra}} associated to $(X,r)$.
\end{definition}
\begin{example}
\label{trivialsolex} For arbitrary nonempty set $X$ with $|X|\geq 2$, \emph{the
trivial solution} $(X, r)$ is defined as $r(x,y)=(y,x),$ for all
$x,y \in X.$ It is clear that $(X,r)$ is the trivial solution
\emph{iff} ${}^xy =y$, and $x^{y} = x,$ for all $x,y \in X,$ or
equivalently $\Lcal_x= \id_X =\Rcal_x$ for all $x\in X.$ In this case
$S(X,r)$ is the free abelian monoid, $G(X,r)$ is the free abelian
group,  $\Acal(k,X,r)$ the algebra of commutative polynomials in
$X$, and $\Gcal(X,r) = \{\id_X\}$.
\end{example}

\begin{remark} \label{ybe} 
Suppose  $(X,r)$ is a nondegenerate quadratic set. 
It is well known, see for example \cite{TSh07},
that $(X,r)$   is a braided set (i.e. $r$ obeys the YBE)
 {\em iff} the following conditions hold
\[
\begin{array}{lclc}
 {\bf l1:}\quad& {}^x{({}^yz)}={}^{{}^xy}{({}^{x^y}{z})},
 \quad\quad\quad
 & {\bf r1:}\quad&
{(x^y)}^z=(x^{{}^yz})^{y^z},
\end{array}\]
 \[ {\rm\bf lr3:} \quad
{({}^xy)}^{({}^{x^y}{(z)})} \ = \ {}^{(x^{{}^yz})}{(y^z)},\]
 for all $x,y,z \in X$.

Clearly, conditions {\bf l1}, imply that for each nondegenerate
braided set $(X,r)$ the assignment: $x \longrightarrow \Lcal_x, x
\in X$ extends canonically to a group homomorphism
\begin{equation}
\label{Lcal} \Lcal:G(X,r) \longrightarrow \Sym(X),
\end{equation}
 which defines  \emph{the canonical left action of $G(X,r)$  on the set  $X$}.
Analogously {\bf r1} implies \emph{a  canonical right action of
$G(X,r)$  on   $X$}.
\end{remark}
\begin{definition} \cite{TSh07}
\label{Gcaldef} Let $(X,r)$ be a nondegenerate braided set,
$\Lcal:G(X,r) \longrightarrow \Sym(X)$ be  the canonical group
homomorphism defined via the left action. The image  $\Lcal(G(X,r))$
is denoted by $\Gcal(X,r).$ We  call it \emph{the (permutation)
group of left actions}.
\end{definition}
The permutation group  $\Gcal(X,r)$ will be of particular importance
in the present paper. 
\begin{remark}
\label{Gcal=LcalS}
Suppose $X$ is a finite set, then $\Gcal=\Lcal(S(X,r))$. Indeed, $\Gcal$ is generated as a semigroup by the images
$\Lcal_x$ of all $x\in X$ and their inverses. Each $\Lcal_x$ is a permutation of finite order, say $m_x$.
Then its inverse (in $\Gcal$) is simply $(\Lcal_x)^{m_x-1} $. So $\Gcal$ is generated as a semigroup by the set 
$\{\Lcal_x, \Lcal_{(x^{m_x-1})} \mid x \in X\} \subseteq \Lcal(S(X,r))$.
\end{remark}

The following conditions were introduced and studied in \cite{T04, TSh07, TSh08}:

\begin{definition}
\label{lri&cl} 
Let $(X,r)$ be a quadratic set.
\begin{enumerate}
\item \cite{T04, TSh07} 
$(X,r)$ is
called {\em cyclic} if the following conditions are satisfied
\[\begin{array}{lclc}
 {\rm\bf cl1:}\quad&  {}^{y^x}x= {}^yx \quad\text{for all}\; x,y \in
 X;
 \quad&{\rm\bf cr1:}\quad &x^{{}^xy}= x^y, \quad\text{for all}\; x,y \in
X;\\
 {\rm\bf cl2:}\quad
  &{}^{{}^xy}x= {}^yx,
\quad\text{for all}\; x,y \in X; \quad & {\rm\bf cr2:}\quad
&x^{y^x}= x^y \quad\text{for all}\; x,y \in X.
\end{array}\]
 We refer to these
conditions as {\em cyclic conditions}.
\item Condition \textbf{lri} is defined as
 \[ \textbf{lri:}
\quad ({}^xy)^x= y={}^x{(y^x)} \;\text{for all} \quad
x,y \in X.\] 
In other words \textbf{lri} holds if and only if
$(X,r)$ is nondegenerate and $\Rcal_x=\Lcal_x^{-1}$ and $\Lcal_x =
\Rcal_x^{-1}$
\end{enumerate}
\end{definition}
In this paper the class of \emph{nondegenerate square-free symmetric sets of finite order} 
will be of special interest.
The following result is extracted from \cite{TSh07}, Theorem 2.34, where more equivalent conditions are 
pointed out. Note that in our considerations below (unless we indicate the contrary) 
the set $X$ is not necessarily  of finite order.
\begin{fact} 
\label{basictheorem} \cite{TSh07}. Suppose $(X,r)$ is
nondegenerate, involutive and square-free quadratic set (not necessarily finite). Then the
following conditions are equivalent:
(i)  $(X,r)$ is a set-theoretic
solution of the Yang-Baxter equation; (ii)  $(X,r)$
satisfies {\bf l1}; 
(iii)  $(X,r)$ satisfies {\bf r1};
(iv)
$(X,r)$ satisfies {\bf lr3}. 

In this case $(X,r)$ is cyclic and satisfies  {\bf lri}.
\end{fact}
\begin{corollary}\label{constructivecor}
Every nondegenerate  square-free symmetric set  $(X,r)$ is uniquely determined by the
left action $\Lcal: X\times X \longrightarrow X,$ more precisely,
\[
r(x,y) = (\Lcal _x(y), \Lcal^{-1}_y(x)).
\]
Furthermore it  is cyclic.
\end{corollary}

Let $(X,r)$ be a braided set.
Clearly, if $Y$ is an $r$-invariant subset  of $(X,r)$, $r$ induces a solution $(Y, r_Y),$
where $r_Y= r_{\mid Y\times Y}$. We call $(Y, r_Y)$ \emph{the restricted solution (on $Y$)}. 
Suppose $(X,r)$ is square-free symmetric set.
Then each $G$-orbit under the left action of $G$ on $X$ is also a right $G$-orbit
and therefore it is an $r$-invariant subset. 

\subsection{Nondegenerate square-free symmetric sets of finite order.}

In the case when  $(X,r)$ is a   square-free  symmetric sets of finite order the
algebras $\Acal(X,r)$ provided new classes of Noetherian rings
\cite{T94,T96}, Gorentstein (Artin-Schelter regular) rings
\cite{T96Preprint,T00,T04s} and so forth. Artin-Schelter regular
rings were introduced in \cite{AS} and are of particular interest.
The algebras $\Acal(X,r)$ are similar in spirit to the quadratic
algebras associated to linear solutions particularly studied in
\cite{Man}, but have their own remarkable properties. The semigroups $S(X,r)$ for a general braided set
$(X,r)$  were studied particularly in \cite{TSh07} with
a systematic theory of `exponentiation' from the set to the
semigroup by means of the `actions' $\Lcal_x,\Rcal_x$ (which in the
process become a matched pair of semigroup actions) somewhat in
analogy with the Lie theoretic exponentiation in \cite{Ma:mat}.

We shall recall some basic facts and recent results needed in this paper.
Suppose $(X,r)$ is a nondegenerate square-free symmetric set of finite order, 
$G= G(X,r), \Gcal=\Gcal(X,r)$ 
in notation as above. 
Note that the group  $G$ acts nontransitively  on $X$. 
This follows from the decomposition theorem of Rump, \cite{Rump}.

Let $X_1, \cdots, X_t$ be the set of orbits of the left action of $G$ on $X$.
These are $r$-invariant subsets of $X$, 
and each $(X_i, r_i), 1 \leq i \leq t$, where $r_i$ is the restriction $r_i= r_{\mid X_i\times X_i}$, is also 
square-free symmetric sets. 

For each $i, 1 \leq i \leq t$, we denote by $\Gcal_i$ the subgroups of $\Sym(X_i)$ 
generated 
by the set of all restrictions
$\Lcal_{x\mid X_i}$:
\[
\Gcal_i = {}_{gr} \langle \Lcal_{x\mid X_i}\mid x \in X \rangle
\]
Note that $\Gcal_i$ acts transitively on $X_i$ for $1 \leq i \leq t.$  
In general $\Gcal_i$ is not a subgroup of $\Gcal$. Rather, there are canonical group surjections $\Gcal\to \Gcal_i$ of restriction to $X_i$ and a canonical inclusion $\Gcal\subseteq \prod_{1\leq i\leq t}\Gcal_i$.

The notions of retraction of symmetric sets and multipermutation solutions were introduced in the general case in \cite{ESS},  where $(X, r)$ is not necessarily finite, or square-free.   In \cite{T04}, \cite{TSh07}, \cite{TSh08} are studied especially the  multipermutation square-free solutions of finite order, we recall some  notions and  results. Let $(X,r)$ be a nondegenerate symmetric set. An
equivalence relation $\sim$ is defined on $X$ as
\[ x \sim y \quad \text{ \emph{iff}} \quad
\Lcal_x = \Lcal_y.\] In this case we also have $\Rcal_x = \Rcal_y,$

We denote by  $[x]$ the equivalence class of $x\in X$, $[X]=
X/_{\sim}$ is the set of equivalence classes.

\begin{lemma}
\label{retractlemma} \cite{TSh07} Let $(X,r)$ be a nondegenerate symmetric set.
\begin{enumerate}
\item
The left and the right actions of $X$ onto itself induce naturally
left and right actions on the retraction $[X],$ via
\[
{}^{[\alpha]}{[x]}:= [{}^{\alpha}{x}]\quad [\alpha]^{[x]}:=
[\alpha^x], \;\text{for all}\; \alpha, x \in X.
\]
\item The new actions (as usual)
define a canonical map $r_{[X]}: [X]\times[X]
\longrightarrow [X]\times[X] $
where $r_{[X]}([x], [y])= ({}^{[x]}{[y]}, [x]^{[y]}).$
\item
$([X], r_{[X]})$ is a nondegenerate symmetric set.
Furthermore,
\item  $(X,r)\; \text{cyclic} \Longrightarrow([X],
r_{[X]})\;\text{cyclic}$.
\item
 $(X,r)\; \text{is}\; {\bf lri} \Longrightarrow([X],
r_{[X]})\;\text{is}\; {\bf lri}.$
\item
$ (X,r)\; \text{square-free} \Longrightarrow ([X], r_{[X]}) \;
\text{square-free}.$
\end{enumerate}
\end{lemma}

\begin{definition}\cite{ESS}
The solution $Ret(X,r)=([X], [r])$ is called the \emph{retraction of
$(X,r)$}.
For all integers $m \geq 1$,  $Ret^m(X,r)$ is defined recursively as
$Ret^m(X,r)= Ret(Ret^{m-1}(X,r)).$

$(X,r)$ is  \emph{a multipermutation solution of level} $m$, 
 if $m$ is
the minimal number (if any), such that $Ret^m(X,r)$ is the trivial
solution on a set of one element. In this case we write $mpl(X,r)=m$. By definition $(X,r)$ is \emph{a multipermutation solution of level}
$0$ \emph{iff} $X$ is a one element set. \end{definition}

The following conjecture was made by the first author in 2004.
\goodbreak
 
\begin{conjecture}  \cite{T04}
\begin{enumerate}
\item
Every nondegenerate square-free symmetric set $(X,r)$ of finite order 
is retractable.
\item
Every nondegenerate square-free symmetric set $(X,r)$ of finite order $n$
is a multipermutation solution, with $mpl(X,r)< n$.
\end{enumerate}
\end{conjecture}

A more recent conjecture states
\begin{conjecture}  \cite{T08ini}
Suppose $(X,r)$ is a nondegenerate square-free multipermutation solution
of finite order $n$. Then $mpl(X) < \log_2 n$.
\end{conjecture}
Evidence for this conjecture and more recent results on 
multipermutation square-free symmetric sets in the general case can be
found in   \cite{TP}.

The following results are of significant importance for our paper and can be deduced from the results in \cite{TSh08}. We give a sketch of the proofs

\begin{proposition}\label{sigprop} Let  $(X,r)$ be a  nondegenerate square-free
symmetric set of finite order, $X_i, 1 \leq i \leq t$ the set of
all $G(X,r)$-orbits in $X$ enumerated so that $X_1, \cdots , X_{t_0}$ is the set of all nontrivial orbits (if any). Then the following are equivalent.
\begin{enumerate}
\item \label{mpl2_1} $(X,r)$ is a multipermutation solution of level 2
\item \label{mpl2_2} $t_0 \geq 1$ and for each $ j, 1 \leq j \leq t_0$,  $x,y \in X_j$ implies $\Lcal_x = \Lcal_y.$
\item \label{mpl2_3} $t_0\geq 1$ and for each $x\in X$ the permutation $\Lcal_x$ is an $r$-automorphism, i.e. $\Gcal(X,r) \subseteq Aut(X,r)$. 
\end{enumerate}
\end{proposition}
\proof If there are no nontrivial orbits then  $mpl(X,r)=1$ since all elements act the same way, i.e. $([X],r_{[X]})$ is the 1-element solution. Assuming $t_0\ge 1$,  $mpl(X,r)=2$ means $mpl([X], r_{[X]})=1$, which means  $[{}^xy]=[y]$ for all $x,y \in X$, i.e.  
\begin{equation} \label{mp2} \Lcal_y=\Lcal_{{}^xy},\quad\forall x,y\in X.\end{equation}
Note that for every pair $x,y \in X$, $x, y$ belong to the same orbit $X_i$ \emph{iff} $y =\Lcal_u (x),$ for some $u \in S(X,r).$ This gives  (\ref{mp2})  $\Longleftrightarrow$  (\ref{mpl2_2}). Meanwhile,  \cite[Lemma 2.7]{TSh08} provides (\ref{mp2})  $\Longleftrightarrow$ (\ref{mpl2_3}). \endproof

\begin{theorem}
\label{significantth}  Let  $(X,r)$ be a  nondegenerate square-free
symmetric set of multipermutation level 2 and finite order, and $X_i$ orbits as in Proposition~\ref{sigprop}.  Let $(X_i, r_i),1 \leq i \leq t $ be the restricted solution. Then:
\begin{enumerate}
 \item  \label{mpl2_4} $\Gcal(X,r)$ is a nontrivial abelian group.
\item \label{mpl2_5} Each $(X_i,r_i), 1 \leq i \leq t_0$ is a trivial solution. Clearly in the case $t_0 < t,$ each $(X_j,r_j)$, with $t_0 \leq j \leq t$ is a one element solution.
\item \label{mpl2_6} For any  ordered pair $i, j, 1\leq i\leq t_0, 1\leq j\leq t,$ such that $X_j$ acts nontrivially on $X_i$, every $x\in X_j$ acts via the same  permutation $\sigma^j_i\in \Sym(X_i)$ which is a product of disjoint cycles of equal length $d=d^j_i$
\[
\sigma^j_i= (x_1\cdots x_d) (y_1\cdots y_d)\cdots (z_1\cdots z_d),
\]
where  each element of $X_i$ occurs exactly once. Here  $d^j_i$  is an invariant of the pair $X_j, X_i$.
\item\label{mpl2_7} $X$ is \emph{a strong twisted union} $X= X_1\natural X_2 \natural \cdots \natural X_t$, (see the definition in \cite{TSh07}).
\end{enumerate}
\end{theorem}

\begin{proof} The proof of part (\ref{mpl2_4}) is given immediately following \cite[Lemma~5.25]{TSh08}. Briefly,  write (\ref{mp2}) as $\Lcal_y=\Lcal_{y^x}$ for all $x,y$ using {\bf lri} (see Fact~\ref{basictheorem}), and then use {\bf l1}. By Proposition~\ref{sigprop}(\ref{mpl2_2}) and $(X,r)$ square-free one has ${}^yx= {}^xx=x$, for all $x,y \in X_i$ which implies  part (\ref{mpl2_5}) here. The implication $mpl(X,r) = 2 \Longrightarrow$ (\ref{mpl2_6})  can be read off from  \cite[Thms. 5.22, 5.24]{TSh08}. The  reader can deduce it also from Proposition~\ref{sigprop}(\ref{mpl2_2}) and the fact that $\Gcal$ is an abelian group of permutations acting transitively on each orbit $X_i$. $mpl(X,r) = 2 \Longleftrightarrow$  (\ref{mpl2_7}) follows from \cite[Thm. 5.24]{TSh08}. \end{proof}

We note that the converse to part (\ref{mpl2_4}) stated in \cite[Thm. 5.24]{TSh08} is incorrect (there was  a gap in the proof in this direction). Indeed, $\Gcal$ a nontrivial abelian group does not  require multipermutation level $2$, see  \cite{TP} for an example.

\section{Representations of finite abelian groups and diagonalization over $\mathbb{C}$}
In this section we shall recall some facts from representation theory of finite abelian groups. 
$\mathbb{C}^{\times}$ will denote the multiplicative group of nonzero complex numbers.
 $C_p$ will denote a cyclic group of finite order $p$, where $p$ is a  positive integer.
We shall write also $C_p=\langle u\rangle$ to denote that $C_p$ is generated by $u$. 
$A$ will denote an abelian group of finite order $N$.  We recall some  facts, from the representation theory of finite abelian groups, which can  be extracted from the literature, 
 see for example  \cite{PeterAlg} and \cite{PeterPerm}.
\begin{theorem}
\label{directsumth}
  Let  $A$ be a finite abelian group. Then 
\begin{equation}
\label{directproduct1}
A = C_{p_1}\times C_{p_2} \times \cdots \times C_{p_s},
\end{equation}
where $p_1 \cdots, p_s$ are positive integers, $p_i \mid  p_{i+1},$ for $1 \leq i \leq s-1.$ 
\end{theorem}
It is not difficult to verify that   the integer $s$ in the presentation (\ref{directsumth})  
is\emph{ minimal}
with the property that $A$ can be presented as a direct product of $s$ cyclic groups.
In addition the condition that $p_i$ divides $p_{i+1},$ for $1 \leq i \leq s-1$ guarantees that the  presentation 
(\ref{directsumth}) is unique, 
so Theorem \ref{directsumth}
is a stronger version of the 
\emph{Basis theorem for abelian groups}, (where originally the condition  $p_i \mid  p_{i+1},$  is not imposed)
see \cite{douglas}.  

We choose a set of generators 
\begin{equation}
\label{B}
\mathbb{B}= \{u_{1},\cdots ,
u_{s}\}
\end{equation}
where $C_{p_i} = \langle u_i\rangle,  1 \leq i \leq s.$  
 Every element $a \in A$ can be presented as
\begin{equation}
\label{coordinates}
a = u_1^{m_1}u_2^{m_2}\cdots u_s^{m_s}, \quad\text{with}\quad 0 \leq m_i \leq p_i-1, \quad \forall i, 1 \leq i \leq s.  
\end{equation}
where 
 the multi-index $s$-tuple ${\bf m}= (m_1,  m_2, \dots, m_s) $ is uniquely determined by $a$ and we refer to it as the {\em coordinates of $a$ with respect to $\mathbb{B}$}. Here  $\mathbb{B}$ is called \emph{a basis for the abelian group $A$}. If we consider $A$ as a $\mathbb{Z}$-module, with operation $+$, then $\mathbb{B}$ behaves as a free-module basis. 

Denote by $\Lambda = \Lambda(A)$ the multi-indexing set 
\begin{equation}
\label{Lambda}
\Lambda = \{{\bf m}= (m_1, m_2, \cdots, m_s ) \mid  0 \leq m_i \leq p_i-1, \quad  1 \leq i \leq s \}
\end{equation}
Clearly,  $\mid \Lambda \mid = \mid A \mid = N$. Note that, in general, the choice of $\mathbb{B}$ is not unique,
but $\Lambda$ is uniquelly determined by the presentation (\ref{directsumth}). 

In the setting of finite abelian grous the notion of a character is  simple and pleasant to work with. As a reference one can use for example  \cite{PeterPerm}. We recall some basic definitions and results. 

\begin{definition} In the context of abelian groups and for our purposes, a  \emph{character} $\chi$ on an abelian group $A$ can be defined simply as a group homomorphism $\chi: A\to \C^\times$, i.e. as a 1-dimensional representation.
\end{definition}
Note that these are {\em irreducible characters} from the point of view of the general theory (and that in general a character on an arbitrary group is defined as the trace of a representation). We will not need the general theory here.

Since $A$ is a finite group of order $N$, then  every $a \in A$, satisfies $a^N = 1,$ 
so $\chi(a^N)= \chi(a)^N = 1,$  therefore $\chi(a)$ is a root of unity of order dividing $N$.

The set $A^{*}$ of (irreducible) characters is a group called the dual of $A$. Suppose $p$ is a positive integer, $C_p= \langle g \rangle $. Let $\theta$ be a primitive $p$-th root of unity.
Then the dual $C_p^{*}$ is itself a  cyclic group of order $p$
generated by the character $\chi^{(p)},$ with $\chi^{(p)}(g)= \theta,$ so $C_p^{*} \cong C_p.$ 
In this case every character $\chi$ is a group homomorphism
\[
C_p \rightarrow \langle \theta \rangle.
\]
More generally, if $A$ is a finite abelian group, presented as (\ref{directsumth}), then
\begin{equation}
\label{A*}
A^{*} =  C_{p_1}^{*}\times C_{p_2}^{*} \times \cdots \times C_{p_s}^{*} \cong C_{p_1}\times C_{p_2} \times \cdots \times C_{p_s} = A.
\end{equation}
Note that the isomorphism $A^{*}\cong A $ is not canonical, since it depends on the decomposition of $A$ as a direct product and on the choice of primitive $p_k$-th roots of unity $\theta_k$ for each $1\le k\le s$. Once this choice is made, we can use our same indexing set $\Lambda$ to label elements of $A^*$:
\begin{equation}\label{chicoord}  \chi(u_k)=\theta_k^{\eta_k}\end{equation}
where $\eta=(\eta_1,\cdots,\eta_s)\in \Lambda$ denotes the $s$-tuple of coordinates of $\chi=\chi_\eta$.

Suppose now that $A\subseteq \Sym(Z)$ is an abelian group of permutations acting \emph{transitively} on a  finite set $Z$. Under the assumptions and notations as above, we  recall some elementary properties which will be useful for our results, see \cite{PeterPerm}, \cite{PeterAlg}. 
\begin{fact}
\begin{enumerate}
\item
For every $x\in Z$ the stabilizer ${\rm Stab}(x)$ is the identity permutation $e=\id_Z.$
\item
Let $x_1 \in Z,$ then the map $A \longrightarrow Z$ defined as $A \ni\pi \mapsto \pi(x_1)$ is a bijection.
Therefore we can consider $Z$ as the set
\[
Z= \{\pi(x_1) \mid \pi \in A \}
\]
In particular,
\[
 \mid Z \mid   = \mid A \mid = N = \prod_{1 \leq k \leq s}p_k
\]
\end{enumerate}
\end{fact}

Fix an arbitrary $x_1 \in Z$. We know that its stabilizer is trivial. An arbitrary  $z \in X,$ can be presented as $z= \pi(x_1),$ for some $\pi \in A,$
since the group $A$ acts transitively on $X$. This permutation $\pi$ is uniquely determined by $z$, indeed 
$\pi (x_1)= \rho (x_1) $  implies $\rho^{-1} \pi (x_1) = x_1,$
and therefore $\rho^{-1} \pi \in {\rm Stab}( x_1) = \{ e \}$. Hence $\rho= \pi.$ Hence there is a 1-1 correspondence between $Z$ and $A$ provided by $z=\pi(x_1)$. Hence the coordinates ${\bf m}={\bf m}_\pi=(m_1,\cdots,m_s)$ of $\pi\in A$ in our basis can also be viewed as coordinates of the corresponding element 
\[ Z=\{x_{\bf m}\ |\ {\bf m}\in\Lambda\};\quad x_{{\bf m}_\pi}=\pi(x_1)\]
In this way the set $\Lambda$ also indexes elements of $Z$.

Finally, we often use $N \times N$ permutation matrices, so it will be convenient to order $\Lambda$ linearly. This induces an enumeration of $A,A^*,Z$ as we have seen.  Let $<$ be the lexicographic order on $\Lambda$. 
$<$ is a linear ordering on a finite set, so we can enumerate  $\Lambda$ accordingly, 
\[
{\bf m}^{(1)} = (0,\cdots, 0) < {\bf m}^{(2)} = (1, \cdots, 0) < \cdots <  {\bf m}^{(N)} = (p_1-1,  \cdots, p_s-1).
\]
We can then refer to our enumerated elements by their position. Thus 
\begin{equation} 
\label{xj}
x_j=x_{{\bf m}^{(j)}}=\pi_j(x_1),\quad \chi_j(u_k)=\theta_k^{m^{(j)}_k}
\end{equation}
where ${\bf m}^{(j)}=(m^{(j)}_1,\cdots,m^{(j)}_s)$.  Clearly, this is equivalent to our previous multi-index
enumeration, and one has  $x_1 < x_2 < \cdots x_N$ as the induced ordering on $Z$.

Let $V$ be the $N$ dimensional $\mathbb{C}$-vector space spanned by $Z.$ 
Then the action of each permutation $\pi \in A$ on $Z$ can be extended
to linear automorphism $T_{\pi}$ of $V$ and the action of $A$ on $Z$ canonically extends
to action of the group algebra $\mathbb{C} A$ on  $V$. 
This induces a matrix representation of $A$ via permutation matrices.

Let $P= P(\pi)$ be the matrix of the automorphism $T_{\pi}$. w.r.t. the basis $Z= \{x_1, x_2, \cdots, x_N\}$.  Then $P$ is a permutation matrix 
with entries
\[
\begin{array}{c}
P_{ij} =\begin{cases} 1 & \text{if}\quad \pi{x_i} = x_j\\ 
                  0 & \text{else} \quad\quad \quad \end{cases} \end{array}
\]
As usual, $\overline{\theta}$    denotes 
the complex conjugate of  $\theta \in \mathbb{C}.$
\begin{theorem}
\label{Diagonalt}
In notation as above, 
 suppose $\chi$ is an irreducible character of $A$. 
Then
\begin{equation}
\label{y-chi}
y_{\chi} = \sum_{1 \leq i \leq N} \overline{\chi(\pi_i)} x_i
\end{equation}
is a simultaneous eigenvector for all $\rho \in A.$
More precisely, there is an equality
\begin{equation}
P(\rho) y_{\chi}  = \chi(\rho)  y_{\chi}. 
\end{equation}
\end{theorem}
\begin{proof}
We know that for all $\pi \in A,$  $\theta= \chi(\pi)  \in \mathbb{C}$ is a root of unity, thus 
$\theta \overline{\theta} = 1,$ and $\overline{\theta}= \theta^{-1}.$ 
This implies 
\begin{equation}
\label{vchi}
 \overline{\chi(\pi)}= (\chi(\pi))^{-1} = \chi(\pi^{-1}).
\end{equation}
Now we use  (\ref{xj}),   (\ref{vchi}) to present  $y_{\chi}$ as
\begin{equation}
\label{vchi1}
y_{\chi} = \sum_{1 \leq i \leq N} \overline{\chi(\pi_i)} x_i= \sum_{1 \leq i \leq N} \chi(\pi_i^{-1})\pi_i(x_1)
\end{equation} 

For each $i, 1 \leq i \leq N$ we set
$\rho \pi_i= \sigma_i$. Then  
\begin{equation}
\label{vchi2}
 \pi_i^{-1}=  \sigma_i^{-1} \rho \quad \chi(\pi_i^{-1})= \chi( \sigma_i^{-1}) \chi(\rho)
\end{equation}
We apply (\ref{vchi}) and (\ref{vchi2}) to deduce   the equalities

\begin{equation}
\begin{array}{clc}
P(\rho)   y_{\chi} &= &P(\rho).\sum_{1 \leq i \leq N} \chi(\pi_i^{-1})\pi_i(x_1) \\
\\
\quad \quad &= & \sum_{1 \leq i \leq N} \chi(\pi_i^{-1})\rho\pi_i(x_1) \\
\\
\quad \quad &= & \sum_{1 \leq i \leq N} [\chi( \sigma_i^{-1}) \chi(\rho)]  \sigma_i(x_1) \\
\\
\quad \quad &= & \chi(\rho) \sum_{1 \leq i \leq N} \chi( \sigma_i^{-1})   \sigma_i(x_1) \\
\\
\quad \quad &= & \chi(\rho) y_{\chi} \quad \quad \quad
\\
\\
\end{array}
\end{equation}
The last equality follows from the equalities of  sets
\[
\{\sigma_i, 1 \leq i \leq N\} = \{\pi_i, 1 \leq i \leq N\}.
\]  
\end{proof}

\section{Complete Datum of multipermutation  solution of level two and explicit formulae for the diagonlization}
\label{SecDatum}
In this section $(X,r)$ denotes a finite nondegenerate
square-free symmetric set of multipermutation level $2$. We shall use the notation of Section 2.
As usual   
$G= G(X,r)$,  $S= S(X,r)$ and  $\Gcal=\Gcal(X,r),$ are respectively the associated YB group, YB monoid, and the permutation group of $(X,r)$.
$X_1, \cdots, X_t $  are the $G$-orbits or equivalently the $\Gcal$-orbits in $X$.

By the Decomposition Theorem 
\cite{Rump}, $t \geq 2.$
By Proposition~\ref{sigprop} there is at least one nontrivial orbit 
$X_i$, 
and we enumerate the orbits so that 
$X_1\cdots X_{t_0}$ is the set of all nontrivial orbits in $X$ ($1\leq t_0\leq t$). 
The restrictions $(X_i, r_i)$ $1\leq i 
\leq t_0$ are trivial solutions and  all elements $x \in X_i$ are
equivalent. In the case when $t_0 < t$, $(X_i, r_i)$ are one element solutions for all
$t_0 < i \leq t.$

We recall from Section~2 that for $1 \leq i \leq t$, $\Gcal_i$ is the subgroup of $\Sym(X_i)$ generated by the set 
$\{\Lcal_{x\mid X_i}\mid x \in X\}.$  As sets we also have 
\[
\Gcal_i=
 \{\Lcal_{u\mid X_i}\mid u \in G(X,r) \}  = \Gcal_{\mid X_i}.\]
 The $\Gcal_i$ are 1-element groups for $t_0< i\leq t$.

\begin{remark}
\label{Gcal=LcalSremark} By Remark \ref{Gcal=LcalS} every element $g$ of $\Gcal$ (respectively of $\Gcal_i$) has
a presentation as $g = \Lcal_u, u \in S$ (respectively as $g = \Lcal_{u\mid X_i}, u \in S$). \end{remark}

Under our assumption the group $\Gcal$ is abelian. Hence each $\Gcal_i$ is an abelian group 
of permutations acting transitively on $X_i$,  $1 \leq i \leq t_0.$ 
It is known, see Section 3,  that in this case there is an equality of orders $\mid\Gcal_i\mid= \mid X_i\mid.$

We shall use the  results of the previous section
(adjusted to our concrete notation for the left actions).
In particular,  for each $i, 1 \leq i \leq t_0,$  $\Gcal_i$ is a direct product of cyclic subgroups (see Theorem \ref{directproduct1})
\begin{equation}
\label{basis11}\Gcal_i= (\Lcal_{u_{i1}\mid X_i})\times
(\Lcal_{u_{i2}\mid
 X_i})\times\cdots\times (\Lcal_{u_{is_i}\mid X_i}),
\end{equation}
where for $1 \leq k \leq s_i,$ the orders $p_{ik}$ of $\Lcal_{u_{ik}\mid X_i}$ satisfy
\begin{equation}
\label{divisibilityoforders1} 
p_{ik} \mid  p_{ik+1}, \quad 1 \leq k \leq s_i-1. 
\end{equation}
Moreover, for each $i, 1 \leq i \leq t_0,$  $s_i,$ is  the least integer such that $\Gcal_i$ is generated
by $s_i$ elements, and by the above remark we can without loss of generality choose all 
$u_{ik}\in S=S(X,r), 1 \leq k \leq s_i$  (but not necessarily
$u_{ij} \in X$.)   The integer $s_i$ and the finite sequence $p_{i1},  \cdots,  p_{is_i}$  in (\ref{basis11}),  
(\ref{divisibilityoforders1}) are invariants of the data $(X_i, r)$ and are uniquelly determined 
for each $i, 1 \leq i \leq t_0$. Note that the set $\{\Lcal_{x\mid X_i}\mid
x\in X\}$ generates the group $\Gcal_i,$ but it is not
necessarily a minimal set of generators.

\begin{definition}
\label{basis2} For $1\leq i \leq t_0$ we fix a  presentation as in
(\ref{basis11}). We call the set $\mathbb{B}_i= \{u_{i1},\cdots ,
u_{is_i}\}$ \emph{a basis of labels for $\Gcal_i$}. We shall assume that $u_{ik}\in S, 1\le k\le s_i$ (see Remark~\ref{Gcal=LcalSremark}).
\end{definition}

We know that each $g\in\Gcal_i$ has a unique presentation in terms of the basis elements:
\begin{equation}
\label{m_w}
g =\prod_{1 \leq k \leq  s_i} (\Lcal_{u_{ik}\mid
X_i})^{m_{g,k}},\quad \text{where}\quad 0 \leq m_{g,k} \leq p_{ik}- 1.
\end{equation}

\begin{definition}
\label{coordinatesdef}
We call the $s_i$-tuple $\bf{m}_g=(m_{g,1}, \cdots, m_{g,s_i})$ the \emph{coordinates of} $g\in \Gcal_i$.
\end{definition}

Consider the vector spaces   $V= {\rm Span}_{\mathbb C} X$, $V_i= {\rm Span}_{\mathbb C} X_i, 1 \leq i \leq t$. Clearly (since $X$ is a basis of $V$), for each $u \in S$ the map $\Lcal_u \in \Sym(X)$
extends canonically to
a  linear automorphism on $V,$ denoted $\Tilde{\Lcal}_u .$
It is a standard fact that
these can be extended to representations of the group algebra
${\mathbb C}\Gcal$ on $V=Span X,$ and since each $V_i, 1 \leq i \leq s$ is invariant,
this representation canonically induces representations of ${\mathbb C}\Gcal_i$ on $V_i.$

We shall apply the results of the previous section
to find a basis 
$Y_i$ of $V_i$ , $1 \leq i \leq t_0$
on which $\Gcal_i$ acts diagonally. In the cases $t_0 < j \leq t$ we simply set $Y_j = X_j.$ 
Then, clearly, the set $Y=\bigcup_{1 \leq i \leq t} Y_i$  is
a basis of $V$ on which the group $\Gcal$ acts diagonally.

To each multipermutation solution  $(X,r)$ with $mpl(X,r)=2$  and in notation as above, we
associated the following  datum $\mathbb{D}$ below. The results Theorems \ref{theoremA}, \ref{maintheorem}, and Corollary \ref{corA}   will give: i) Explicit coordinates of a new basis $Y_i$ w.r.t initial basis $X_i$ of $\C X_i$ in terms of the datum $\mathbb{D}$. The new basis will have the feature that the linear solution $R$ associated to $r$ is diagonal. ii) Explicitly presented form of the diagonal coefficients  of $R$.

\begin{definition}\label{datum} For $1 \leq i \leq t_0$, we define \emph{the complete datum $\mathbb{D}_i$ of $X_i$} as $(X_i,\mathbb{B}_i,\Pcal_i,\Theta^i,\Mcal_i)$ where:
\begin{enumerate}
\item $\mathbb{B}_i= \{u_{i1},\cdots , u_{is_i}\}$  is a basis as in Definition
\ref{basis2}. 
\item $\Pcal_i= \{p_{i1},\cdots ,p_{is_i}\}$, where $p_{ik}$ is the
order of $\Lcal_{u_{ik}\mid X_i}, 1 \leq k\leq s_i$. Note that
\[ N_i=p_{i1}\cdots p_{is_i}=\mid X_i\mid =\mid \Gcal_i\mid\]
by the results of the previous section. 
\item $\Theta ^i= \{\theta_{i1},\cdots ,\theta_{is_i}\} $ where $\theta_i$  is a chosen primitive $p_{is_i}$--th root of unity
and 
\[
\theta_{ik}= (\theta _i)^{\frac{p_{is_i}}{p_{ik}}}, \quad 1 \leq k \leq s_i.
\]
Clearly,  each $\theta_{ik}$ is a primitive $p_{ik}$--th root of unity. 
\item $\Mcal_i=\{{\bf m}^j_i\ |\ 1\le j\le t\}$ where ${\bf m}^j_i=(m^j_{i1},\cdots,m^j_{is_i})$ denotes the coordinates in our basis of the element $\sigma^j_i\in \Gcal_i$ for the action of any element of $X_j$ on $X_i$ (see Theorem~\ref{significantth}). 
\item $\mathbb{D}=\mathbb{D}(X,r)=(\mathbb{D}_1, \cdots, \mathbb{D}_t)$ is \emph{a complete datum of} $(X,r).$
\end{enumerate}
\end{definition}

In the particular case, when $i>t_0$, i.e. $X_i$ is a one element set  and   (for completeness) we define
$\mathbb{D}_i =(X_i, \{1\}, \{1\}, \{1\}, \{1\} )$.


 From the $\Pcal_i$ part of the datum for each $i$, $1\leq i\leq t$, we also have an associated space
\[ \Lambda_i= \{{\bf m}=(m_{1},\cdots ,m_{s_i})\mid
0 \leq m_{k} \leq p_{ik}- 1,\ 1 \leq k \leq s_i\}\]
for the range of coordinates in our basis $\mathbb{B}_i$.  Clearly, $\mid \Lambda_i \mid =N_i$. Every element  of $\Gcal_i$ can be written with coordinates ${\bf m}_g\in \Lambda_i$ as in Definition~\ref{coordinatesdef}. As explained in Section~3 this induces also coordinates on $X_i$ once we fix an element  $x_{i1}\in X_i$, and it induces coordinates for $\Gcal_i^*$ given the data $\Theta^i$. Thus, given $\eta=(\eta_1,\cdots,\eta_{s_i})\in \Lambda_i$ the corresponding character is  
\begin{equation}
\label{chi-eta}
\chi(\Lcal_{u_{i1}\mid X_i})=\theta_{i1}^{\eta_1},  \quad \chi(\Lcal_{u_{i2}\mid X_i})=\theta_{i2}^{\eta_2}, \quad \cdots,\quad \chi(\Lcal_{u_{is_i}\mid X_i})=\theta_{is_i}^{\eta_{s_i}}\end{equation}

\begin{notation}
Let $\eta\in \Lambda_i,$ and $\chi_{\eta}\in \Gcal_i^{*}$  the corresponding character as in \ref{chi-eta}. The simultaneous eigenvector vector $ y_{\chi_{\eta}}$, see Theorem~\ref{Diagonalt} and (\ref{y-chi}), 
will be denoted  $y_{\eta}^i$.
\end{notation}

Finally, in order to have explicit formulae, we enumerate $\Lambda_i$ as in Section~3. We enumerate the corresponding elements of $\Gcal_i$ as $\id=\pi_1<\cdots<\pi_{N_i}$ and the elements of $X_i$ as $x_{i1}<x_{i2}<\cdots<x_{iN_i}$.  As in the proof of Theorem   \ref{Diagonalt}   we then have explicitly 
\begin{equation}\label{yetachi}
y_{\eta}^i = y_{\chi_{\eta}}=\sum_{1 \leq j \leq N_i} \chi_{\eta}(\pi_j^{-1}) x_{ij}
\end{equation}
for the explicit change from a basis $X_i$ to the new basis $Y_i=\{y^i_\eta\ |\ \eta\in \Lambda_i\}$.  Here $x_{ij}=\pi_j(x_{i1})$. We can also write this formula as a sum over $\pi\in\Gcal_i$ (without enumeration). 

\begin{lemma}
\label{y_eta}
 In notation and enumeration as above, let $\eta = (\eta_1, \cdots, \eta_{s_i}) \in \Lambda_i.$ 
Denote 
\[
\lambda_1 = (\theta_{i1})^{\eta_1}, \quad \lambda_2 = (\theta_{i2})^{\eta_2}, \cdots, 
\quad \lambda_{s_i} = (\theta_{is_i})^{\eta_{s_i}}. 
\]
The simultaneous eigenvector $y_{\eta}^i$ 
has the following explicit coordinates with respect to the basis $X_i$ (given in the form  of a vertical "block-vector"):
\begin{equation}
\label{beta} y^i_{\eta}=\left[
\begin{array}{c}
\\ 
(\lambda_{s_i})^{p_{is_i}-1}B_{s_i-1}
\\
\\ 
\hline
\\ 
\vdots
\\
\\
 \hline
\\
(\lambda_{s_i})^2B_{s_i-1}
\\
\\  
\hline
\\
(\lambda_{s_i})B_{s_i-1}
\\
\\ 
 \hline
\\
B_{s_i-1}
\\
 \end{array}\right],
\end{equation}
where the vectors $B_k$ are determined recursively as follows.
\[
 B_1=\left[\begin{array}{c} \lambda_1^{p_{i1}-1}
\\
\\
\vdots
\\
\\
\lambda_1^2\\
\\
\lambda_1\\
\\
1
 \end{array}\right],
\]
and for $2\leq k \leq s_i,$ one has
\[
B_k=\left[\begin{array}{c}
\\
(\lambda_k)^{p_{ik}-1}B_{k-1}
 \\
\\  
\hline
\\ 
\vdots
\\
\\
 \hline
\\
(\lambda_k)^2B_{k-1}
 \\ 
\\ 
\hline
\\
\lambda_kB_{k-1}
\\
\\  
\hline
\\
B_{k-1}
\\ 
 \end{array}\right]
\]
\end{lemma}
In particular, for $\eta^1= (0, 0, \cdots , 0)$ one has
\begin{equation*}
 y^i_{\eta^1}= \left[\begin{array}{c} 1
\\
1
\\
\vdots
\\
1\\
1
 \end{array}\right],
\end{equation*}
The following theorem is mainly an elaboration of  Theorem \ref{Diagonalt} in terms of 
the fixed datum and the enumeration of $X$ induced from the action of $\Gcal$.

\begin{theorem}
\label{theoremA} Let $(X,r)$ be a finite  nondegenerate
square-free symmetric set and notations as in Section~2.  We assume that $(X,r)$ is a nontrivial solution.
Let $V=\C X=\oplus_i V_i$ where $V_i=\C X_i$. Then
\begin{enumerate}
\item
\label{TA1}
The left action of of the finite group $\Gcal$ on $X$,
extends canonically  to a representation of the group algebra
${\mathbb C}\Gcal$ on $V$. 
Each $V_i, 1 \leq i \leq t,$ is invariant under this representation. Moreover, suppose there exists a basis $Y$ of $V$ on which $\Gcal$ acts diagonally. Then $\Gcal$ is abelian.
\item Suppose that $(X,r)$ has $mpl(X,r)=2$ and datum $\mathbb{D}$.   Then: \begin{enumerate}
\item
\label{TA3}
The set $Y_i=\{y^i_\eta\ |\ \eta\in \Lambda_i\}$ is a basis of $V_i$ on
which $\Gcal_i$ acts diagonally.
The set $Y=\bigcup Y_i,$ is
a basis of $V$ on which the group $\Gcal$ acts diagonally.
\item
\label{TA4}
 The elements of $Y_i$ are simultaneous eigenvectors of  all
$g \in \Gcal_i,$ with eigenvalues $\mu^g_{i,\eta}$ determined from the
coordinates ${\bf m}_g$,
\begin{equation}
\label{theq41} g. y_{\eta}^i=\mu_{i,\eta}^{g}y_{\eta}^i,
\quad \mu_{i,\eta}^{g}=\chi_\eta(g)=\prod_{1 \leq k\leq
s_i}(\theta_{ik}^{\eta_k})^{m_{g,k}} .
\end{equation}
\end{enumerate}\end{enumerate}
\end{theorem}
\begin{proof}
Part (\ref{TA1}) collects some classical facts.  Note that by its definition as a subgroup of $\Sym(X)$, the group $\Gcal$ acts faithfully on $V$. In the basis $Y$ the matrices for the elements of $\Gcal$, being diagonal, all commute. Since they also form a faithful representation, the group $\Gcal$ must be abelian.  Part (2) is an elaboration of Theorem \ref{Diagonalt} in the present context.\end{proof}

The next corollary follows straightforwardly from Theorem \ref{theoremA}.
Using direct information from the datum $\mathbb{D}$ it gives  explicitly  the coefficients which occur in the Main Theorem 
\ref{maintheorem}. 

Recall that each ordered pair $X_j, X_i$ defines uniquely the $s_i$-tuple ${\bf m}_i^j =(m_{i,1}^j, \cdots ,m_{i,s_i}^j) \in \Mcal_i$. (These are the coordinates of any element $\Lcal_{x\mid X_i }$ when $x \in X_j$).   
\begin{corollary} 
\label{corA} In the hypothesis of Theorem \ref{theoremA}.
Let $i,j$ be  integers $1 \leq i,j\leq t$.  Then all $x \in X_j$ act in the same way as
\begin{equation}
\label{theq1} \Lcal_x y_{\eta}^i=\mu_{i, \eta}^{j}y_{\eta}^i;\quad \mu_{i, \eta}^{j}=\chi_\eta(\sigma^j_i)= \prod_{1 \leq k\leq
s_i}(\theta_{ik}^{\eta_k})^{m_{i,k}^j}.
\end{equation}
where ${\bf m}^j_i\in\Mcal_i$ is from our datum. In particular, 
$\mu_{i, \eta}^{j}=1$, whenever $j=i,$ so
\[
\Lcal_x y_{\eta}^i = y_{\eta}^i,\quad\forall x\in X_i.
\]
\end{corollary}

We are now ready to prove our main theorem. 

\begin{theorem}
\label{maintheorem} Let $(X,r)$ be a finite  nondegenerate
square-free symmetric set of order n,  with $mpl(X,r)=2,$ and datum $\mathbb{D}.$
For each $i, 1 \leq i \leq t,$ let  $Y_i$ be the basis of $V_i$ determined by Theorem \ref{theoremA}.
Then the linear extension $R: V \otimes V
\longrightarrow V \otimes V,$ of the solution $r$ to $kX$ has the form
\[
R(y_{\eta}^i \tens y_{\zeta}^j) = \frac{\mu_{j, \zeta}^{i}}{\mu_{i, \eta}^{j}}\;
y_{\eta}^i \tens y_{\zeta}^j.
\]
\end{theorem}

\begin{proof}
Note first that as a linear  extension of the map $r$, the automorphism $R$ satisfies
\begin{equation}
\label{mtheqR}
R(z \tens x) = {}^zx \tens z^x.
\end{equation}
Note also that since {\bf lri} holds, action from the right is inverse to action from the left, and that under our assumptions this left action has the form ${}^zx=\sigma^j_i(x)$ if $z\in X_j$ and $x\in X_i$. Hence we merely need to compute
\begin{eqnarray*}
R(y^i_\eta\tens y^j_\zeta)&=&\sum_{\pi\in \Gcal_i\atop \rho\in\Gcal_j}\chi_\eta(\pi^{-1})\chi_\zeta(\rho^{-1})R(\pi(x_{i1})\tens \rho(x_{j1}))\\
&=&\sum_{\pi\in \Gcal_i\atop \rho\in\Gcal_j}\chi_\eta(\pi^{-1})\chi_\zeta(\rho^{-1})\ {}^{\pi(x_{i1})}\rho(x_{j1})\tens \pi(x_{i1})^{\rho(x_{j1})}\\
&=&\sum_{\pi\in \Gcal_i\atop \rho\in\Gcal_j}\chi_\eta(\pi^{-1})\chi_\zeta(\rho^{-1})(\sigma^i_j\rho)(x_{j1})\tens ((\sigma^j_i)^{-1}\pi)(x_{i1})\\
&=& \sum_{\pi'\in \Gcal_i\atop \rho'\in\Gcal_j}\chi_\eta(\pi'{}^{-1}(\sigma^j_i)^{-1})\chi_\zeta(\rho'{}^{-1}\sigma^i_j)\rho'(x_{j1})\tens \pi'(x_{i1})\\
&=&\chi_\eta((\sigma^j_i)^{-1})\chi_\zeta(\sigma^i_j) \sum_{\pi'\in \Gcal_i\atop \rho'\in\Gcal_j}\chi_\zeta(\rho'{}^{-1})\rho'(x_{j1})\tens  \chi_\eta(\pi'{}^{-1}) \pi'(x_{i1}) \\ 
&=&{\chi_\zeta(\sigma^i_j)\over\chi_\eta(\sigma^j_i)}y^j_\zeta\tens y^i_\eta\end{eqnarray*}
where we use the definition (\ref{yetachi}) and linearity of $R$. We change the summation over $\pi,\rho$ to summation over $\pi'= (\sigma^j_i)^{-1}\pi$ and $\rho'=\sigma^i_j\rho$ (if one runs over $\Gcal_i$ and $\Gcal_j$ respectively then so does the other). We then use that the characters are group homomorphisms, identifying the result as stated. \end{proof} 

We do not know if only $mpl(X,r)=2$ solutions are diagonalisable in this way but in the next section we give an example in support of that conjecture as well as `forward' examples of Theorem~\ref{maintheorem}. 

\section{Worked out examples }

In this section we illustrate the results on Section~4 on some nontrivial examples. The first example was suggested by Peter Cameron. 

\begin{example}
\label{useofbasisthm}
Let $(X,r)$ be the square-free symmetric set defined as follows.
\[
\begin{array}{lcl}
X &=& \{x_1, x_2, \cdots, x_{12}, a, b\}\\
& &\quad \\
\Lcal_{a}&=&(x_1x_3x_5x_7x_9x_{11})(x_2x_4x_6x_8x_{10}x_{12})
\\
& &\quad \\
\Lcal_{b}&=&(x_1x_4x_7x_{10})(x_2x_5x_8x_{11})(x_3x_6x_9x_{12})
\\
& &\quad \\
\Lcal_{x_i}&=&\id_X \quad 1 \leq i \leq 12.
\end{array}
\]
The graph for this example is shown in Figure~1.
\end{example}
\begin{figure}
\[ \includegraphics[scale=.9]{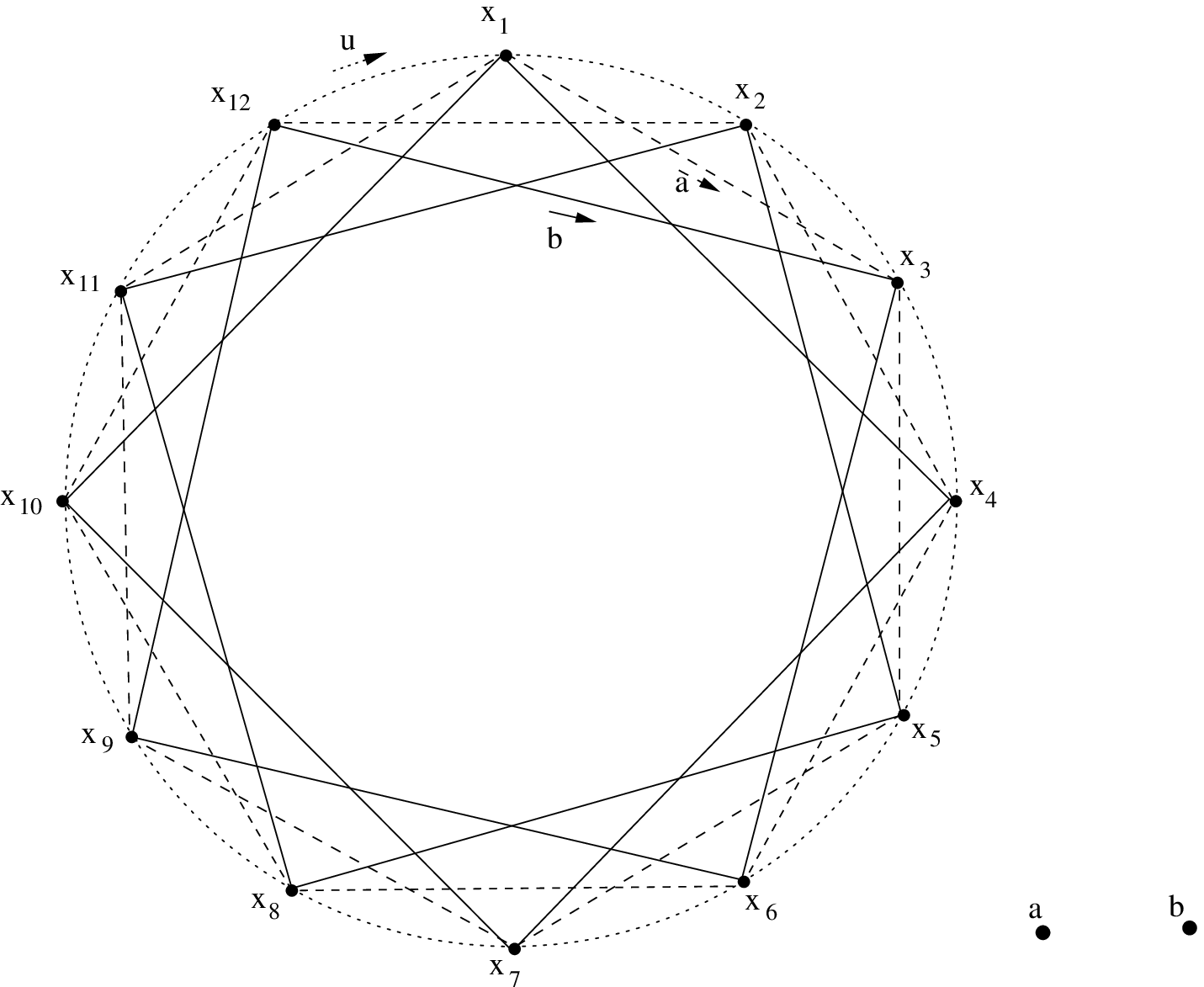}\]
\caption{Graph for Example~\ref{useofbasisthm}. Arrows show left actions on the set $X=\{x_1,\cdots x_{12}\}\cup\{a\}\cup\{b\}$.}
\end{figure}

In this case $\Gcal= {}_{gr}\langle\Lcal_{a},\Lcal_{b} \rangle$, and $X$
splits into $3$ $\Gcal$-orbits:
$X_1=\{x_i \mid 1 \leq i \leq 12\}$, $X_2 = \{a\}, X_3 = \{b\}.$ Clearly, $\Gcal_1 = \Gcal.$
Note that $(\Lcal_a)^3= (\Lcal_b)^2,$ so we can not present $\Gcal_1$ as
a direct product of the cyclic groups
$\langle\Lcal_a\rangle $ and $\langle\Lcal_a\rangle $. One can use as a
basis
for $\Gcal_1$ the element $u= a^2b^3$. More precisely,
\[
\Gcal_1= \langle \Lcal_{u} \rangle, \quad  \Lcal_a =(\Lcal_{u})^{2},
\quad \Lcal_b=(\Lcal_{u})^{3},
\]
$\Lcal_{u}= (x_1 x_2 x_3 x_4 x_5 x_6 x_7 x_{8} x_9 x_{10} x_{11}
x_{12})$ is a cycle of length 12.
The datum $\mathbb{D}(X,r)$ is:
\[
\begin{array}{c}
s_1=1, \quad \mathbb{B}_1= \{u\}, \quad \Pcal_1= \{12\} \\
\\
\Theta ^1= \{\theta\}\quad\text{where $ \theta$ is a chosen primitive
12--th root of unity}
\\
\\
\Mcal_1=\{{\bf m}^2_1,{\bf m}^3_1 \}\quad\text{where} \quad {\bf
m}^2_1=(2), {\bf m}^3_1=(3)
\\
\\
\mathbb{D}_1 = (X_i,\mathbb{B}_1,\Pcal_1,\Theta^1,\Mcal_1),\  \mathbb{D}_2
= (X_2,\{1\},\{1\},\{1\},\{1\}),\  \mathbb{D}_3 =
(X_3,\{1\},\{1\},\{1\},\{1\})
\\
\\
\mathbb{D}=\mathbb{D}(X,r)=(\mathbb{D}_1,\mathbb{D}_2, \mathbb{D}_3)
\\
\\
\Lambda_1 = \{k\mid 0 \leq k \leq 11\}
\end{array}
\]
Here instead of $y_{\eta}^1,\eta \in \Lambda_1 $ we shall simply write
$y_k, 1 \leq k \leq 11,$
$y_1^2= a, y_1^3=b.$
Note that this enumeration agrees with our convention in Section 4. 
Then    $y_k,$  in coordinates with respect to the basis $X_1$ of $V_1$  and written as the  transpose of a row vector is
\[
\begin{array}{c}
y_{k} = \left[ \begin{array}{clclcl} (\theta^k)^{11}& \; \cdots &\;
(\theta^k)^2 &\;\theta^k &\; 1 \\
                 \end{array}\right]^t
 \end{array}
\]
Furthermore, one has
\[
R(a \tens y_k) = \theta^{2k}\; y_k \tens a,  \quad R(b \tens y_k) =
\theta^{3k}\; y_k \tens b \]
\[ R(y_j \tens y_k) = y_k \tens y_j ,\quad \quad R(a \tens b) = b \tens a,\quad
\forall k,j, 1 \leq k,j \leq 12
\]

\begin{example}
\label{moreinterestingex} 
Let 
$(X,r)$ be the square-free symmetric set defined as follows.
\[
\begin{array}{cl}
X = X_1 \bigcup X_2\bigcup X_3, &\quad X_1 = \{a_i\mid 1 \leq i \leq 18\}\\
\\
X_2 = \{b_j\mid 1 \leq j \leq 16\}, &\quad X_3 = \{c_k\mid 1 \leq k \leq 8\}
\end{array}
\]
The solution $r$ is defined via the left actions:
\[
\Lcal_{a_i}=\tau.\upsilon, \quad 1 \leq i\leq 18;\quad 
\Lcal_{b_j}=\rho.\omega, \quad 1 \leq j\leq 16;\quad
\Lcal_{c_k}=\pi.\sigma, \quad 1 \leq k\leq 8 
\]
where the permutations $\pi, \rho, \sigma, \tau, \upsilon, \omega$ are given below.
\[
\begin{array}{clc}
\pi = \pi_0 \pi_1 \pi_2\pi_3\pi_4\pi_5,&\quad\text{where}\quad \pi_k=(a_{3k+1}\; a_{3k+2}\; a_{3k+3}),&\quad 0\leq k \leq 5\\ 
\\
\rho = \rho_1 \rho_2 \rho_3, &\quad\text{where}\quad \rho_m=(a_{m}\; a_{m+3}\; a_{m+2.3}\;\cdots\; a_{m+5.3}),&\quad 1\leq m \leq 3\\ 
\\
\tau = \tau_0 \tau_1 \tau_2 \tau_3, &\quad\text{where}\quad \tau_j=(b_{4k+1}\; b_{4k+2}\; b_{4k+3} \;b_{4k+4}),&\quad 0\leq k \leq 3
\\
\\
\sigma = \sigma_1 \sigma_2 \sigma_3\sigma_4,&\quad\text{where}\quad \sigma_m=(b_{m}\; b_{m+4} \;b_{m+2.4} \;b_{m+3.4}),&\quad 1\leq m \leq 4\\ 
\\
\upsilon= \upsilon_0 \upsilon_1 \upsilon_2\upsilon_3,&\quad\text{where}\quad \upsilon_m=(c_{2k+1} \;c_{2k+2}),&\quad 0\leq k \leq 3
\\ 
\\
\omega = \omega _1 \omega _2, &\quad\text{where}\quad \omega _m=(c_{m}\; c_{m+2}\; c_{m+2.2}\;c_{m+3.2}),&\quad 1\leq m \leq 2
\end{array}
\] The graph for this example is given in Figure~2.
\end{example}
\begin{figure}
\[ \includegraphics[scale=.95]{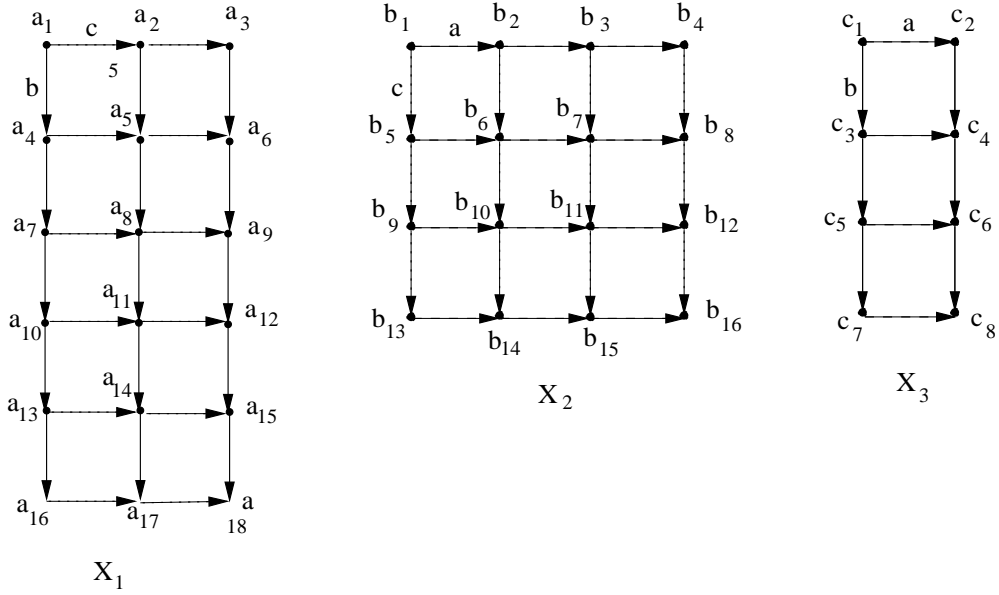}\]
\caption{Graph for Example~\ref{moreinterestingex}. Arrows show left actions on the set $X=X_1\cup X_2\cup X_3$.}
\end{figure}
Clearly, $\Gcal={}_{gr}\langle \Lcal_{a_1}, \Lcal_{b_1}, \Lcal_{c_1}\rangle$  and the $\Gcal$-orbits in $X$
are exactly $X_1, X_2, X_3.$ This time each of the restricted groups $\Gcal_i, 1 \leq\i \leq 3$ is nontrivial and differs from $\Gcal.$ More precisely, one has 
\[
\begin{array}{c}
\Lcal_{a\mid X_1}= e,\quad \Lcal_{b\mid X_1}= \rho,\quad \Lcal_{c\mid X_1}= \pi,\quad\forall a \in X_1, \;b \in X_2,\; c \in X_3 \\
\\
\Gcal_1= \langle \pi \rangle \times  \langle \rho \rangle, \quad
\pi^3=e,\quad \rho^6=e .
\end{array}
\]
In the usual notation one has: 
\[
\begin{array}{c}
s_1=2, \quad \mathbb{B}_1= \{u_{11}= c_1, u_{12}= b_1 \}, \quad \Pcal_1= \{3, 6\}, \\
\\
\Theta ^1= \{\theta_{11}= (\theta)^2, \theta_{12}= \theta \},\quad\text{where $\theta$ is a chosen primitive 6--th root of unity},
\\
\\
\Mcal_1=\{{\bf m}^2_1=(0,1),\quad {\bf m}^3_1=(1,0) \},
\quad
\mathbb{D}_1 = (X_1,\mathbb{B}_1,\Pcal_1,\Theta^1,\Mcal_1); \\
\\
\Lambda_1 = \{(\eta_1, \eta_2) \mid 0 \leq \eta_1 \leq 2,\quad 0 \leq\eta_2\leq 5 \}.
\end{array}
\]
On $X_2$ one has  
\[
\begin{array}{c}
\Lcal_{a\mid X_2}= \tau,  \quad \Lcal_{b\mid X_2}= e, \quad \Lcal_{c\mid X_2}= \sigma, \quad \forall a \in X_1, \;b \in X_2, \;c\in X_3 \\
\\
\Gcal_2= \langle \tau \rangle \times  \langle \sigma \rangle \quad
\tau^4=e\quad \sigma^4 =e 
\end{array}
\]
\[
\begin{array}{c}
s_2=2 \quad \mathbb{B}_2= \{u_{21}= a_1, u_{22}= c_1 \}, \quad \Pcal_2= \{4, 4\} \\
\\
\Theta ^2= \{\theta_{21}= i,\; \theta_{22}= i \}\quad\text{(here we chose as a primitive 4--th root of unity simply $\; i \in \mathbb{C}$}) 
\\
\\
\Mcal_2=\{{\bf m}^1_2=(1,0),\quad {\bf m}^3_2=(0,1) \} \quad \mathbb{D}_2 = (X_1,\mathbb{B}_2,\Pcal_2,\Theta^2,\Mcal_2); \\
\\
\Lambda_2 = \{(\eta_1, \eta_2) \mid 0 \leq \eta_1, \eta_2 \leq 3 \}
\end{array}
\]
Finally,
\[
\begin{array}{c}
\Lcal_{a\mid X_3}= \upsilon,  \quad \Lcal_{b\mid X_3}= \omega,\quad \Lcal_{c\mid X_3}=e, \quad \forall a \in X_1,\; b \in X_2,\; c\in X_3 \\
\\
\Gcal_3= \langle \upsilon \rangle \times  \langle \omega \rangle, \quad
 \upsilon^2=e,\quad \omega^4 =e 
\end{array}
\]
Furthermore,  
\[
\begin{array}{c}
s_3=2, \quad \mathbb{B}_3= \{u_{31}= a_1, u_{32}= b_1 \}, \quad \Pcal_3= \{2, 4\} \\
\\
\Theta ^3= \{\theta_{31}= i^2= -1, \theta_{32}= i \}
\\
\\
\Mcal_3=\{{\bf m}^1_3=(1,0),\quad {\bf m}^3_2=(0,1) \}
\\
\\
\mathbb{D}_3 = (X_3,\mathbb{B}_3,\Pcal_3,\Theta^3,\Mcal_3); \\
\\
\Lambda_3 = \{(\eta_1, \eta_2) \mid 0 \leq \eta_1\leq 1,  0 \leq  \eta_2 \leq 3 \}
\end{array}
\]
\[
\mathbb{D}=\mathbb{D}(X,r)=(\mathbb{D}_1,\mathbb{D}_2, \mathbb{D}_3)
\]
We can compute all $y_\eta^i$, $1 \leq i \leq 3, \eta \in \Lambda_i,$ via the recursive formulae in Lemma \ref{y_eta}. For example, the 18 entries of  $y_{(1,1)}^1$ written as a transposed row vector are
\[
\begin{array}{cl}
y_{(1,1)}^1& \quad\\
&=\left[ \begin{array}{clclclclclclclclcl} \theta^5\theta^4& \; \theta^5\theta^2& \;
\theta^5& \;\theta^4\theta^4& \; \theta^4\theta^2& \;
\theta^4& \;\theta^3\theta^4& \; \theta^3\theta^2& \;
\theta^3& \;\theta^2\theta^4& \; \theta^2\theta^2& \;
\theta^2& \;\theta\theta^4& \; \theta\theta^2& \;
\theta& \;\theta^4& \; \theta^2& \;
1 \\
\end{array}\right]^t\\ \\
&= \left[ \begin{array}{clclclclclclclclcl} \theta^3& \; \theta& \;
\theta^5& \;\theta^2& \; 1& \;
\theta^4& \;\theta& \; \theta^5& \;
\theta^3& \;1& \; \theta^4& \;
\theta^2& \;\theta^5& \; \theta^3& \;
\theta& \;\theta^4& \; \theta^2& \;
1 \\
\end{array}\right]^t
\end{array}
\]

The above two examples illustrate our main Theorem~\ref{maintheorem} as well as our explicit methods developed in Section~4. In the converse direction to Theorem~\ref{maintheorem}, we can ask:

\begin{question}\label{conjR}
Suppose that $(X,r)$ is a square-free nondegenerate symmetric set of finite
order $n$, $X_i$ orbits of $G(X,r)$, $V=\C X$ and $V_i=\C X_i$.  Suppose that  there exists a basis $Y$ of $V$ with $Y=\cup Y_i$, each $Y_i$ a basis of $V_i$, such that  the linear extension $R$ of $r$ has the form
\[R(y_i\tens y_j) = \mu_{ij} y_j\tens y_i, 1 \leq i,j \leq n ,\quad\mu_{ij}\in\C^\times.\] 
Can we conclude that $(X,r)$ is a multipermutation solution with  $mpl(X,r) =2$?
\end{question}

We provide an example in support of an affirmative answer here. This will  be a square-free symmetric set $(X,r)$ with $mpl(X,r)=3$ for which it is impossible to find a new basis $Y$  of the diagonal form stated. From the list of solutions of order $\leq 4$ in \cite{T94} it is not difficult to verify that $|X|\le 4$ implies that $mpl(X,r)\le 2$, hence such an example must have order at least $|X|=5$, see also \cite{TP}. A method for constructing higher permutation level examples is in \cite{TP} and we use one of these.

\begin{example} \cite{TP}
\label{counterex}
Let $(X,r)$ be the square-free symmetric set of order $5$ defined as
follows.
\[ X = \{x_1, x_2,x_3, x_4, a\} \]
\[ \Lcal_{a}=(x_1x_2x_3x_4), \quad \Lcal_{x_1}=\Lcal_{x_3} = (x_2x_4), \quad \Lcal_{x_2}=\Lcal_{x_4} = (x_1x_3)\]
\[
\Gcal = {}_{gr}\langle \Lcal_{a}, \Lcal_{x_1}, \Lcal_{x_2}\rangle = {}_{gr}\langle \Lcal_{a}, \Lcal_{x_2}\rangle={}_{gr}
\langle (x_1 x_2 x_3 x_4),(x_1 x_3) \rangle \simeq D_8,
\]
where $D_8$ is the Dihedral group of order 8 with generators of order $4$ and $2$ as exhibited.
\end{example}
In this case the orbits are $X_1 =\{x_1,x_2, x_3, x_4\}$,  $X_2 =
\{a\}$.  It is straightforward to see that the retract $Ret(X,r) = ([X],
r_{[X]})$
satisfies 
\[ [X] = \{[x_1], [x_2], [a]\},\quad \Lcal_{[a]} = ([x_1][x_2]), \quad \Lcal_{[x_1]}= \Lcal_{[x_2]} =
\id_{[X]}\]
so  $mpl(Ret(X,r)) = 2$ and thus  $mpl(X,r) = 3$. The graph for this example is shown in Figure~3.
\begin{figure}
\[ \includegraphics{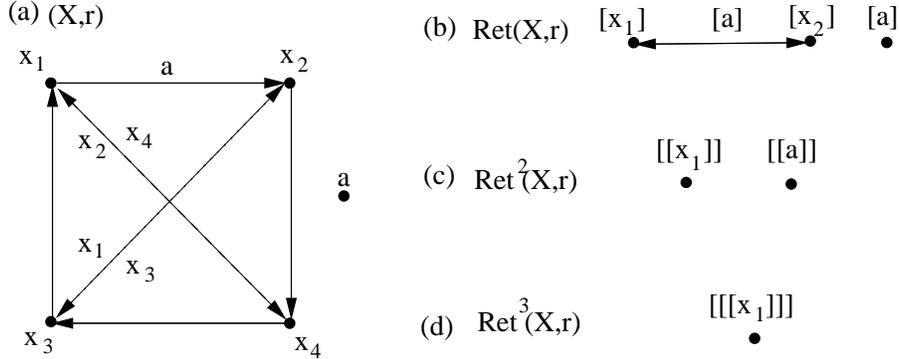}\]
\caption{Graphs for Example~\ref{counterex}. Arrows show left actions on the set $X=\{x_1,x_2,x_3,x_4\}\cup\{a\}$ in part (a). The first retract is shown in part (b). The second retract (c) is the trivial solution. The third retract (d) is the 1-element solution. Hence $mpl(X,r)=3$.}
\end{figure}

Now assume that there is a basis $Y=Y_1 \bigcup \{a\} $  fulfilling the conditions in Conjecture~\ref{conjR}.  Let  $Y_1 = \{y_1, y_2, y_3, y_4 \} \subset V_1$.  So we assume that 
\begin{equation}
\label{R(ay)}
R(a\tens y_i) = \lambda_i y_i\tens a ; \quad \lambda_i
\in \mathbb{C}^{\times},\quad  1 \leq i \leq 4
\end{equation}
and 
\begin{equation}
\label{R(ay)}
R(y_i\tens y_j) = \mu_{ij} y_j\tens y_i ; \quad
\mu_{ij} \in \mathbb{C}^{\times},\quad 1 \leq i \leq 4.
\end{equation}
Each  $y_i\in Y_1$ is some linear combination $y_i=\sum_j \alpha_j x_j$ (say) but the right action on $a$ is trivial, so
\[ R(a\tens y_i)=R(a\tens\sum\alpha_jx_j)=\sum\alpha_j R(a\tens x_j)=\sum \alpha_j\, {}^a(x_j)\tens a=\Lcal_a(y_i)\tens a.\]
Hence we require $\Lcal_a y_i\tens a=\lambda_i y_i\tens a$. Hence $\Lcal_a y_i=\lambda_i y_i$. So each $y_i\in Y$ is an eigenvector of $P(a)$ -- the permutation matrix of the cycle $(x_1 x_2 x_3 x_4)$. Since the latter has  four eigenvectors and these have distinct eigenvalues, this determines $Y_1$ up to normalisation and order. Without loss of generality we can therefore take
\[ Y_1=\left\{y_0=\begin{pmatrix}1\\ 1\\ 1\\1\end{pmatrix},\quad y_1=\begin{pmatrix}-\imath\\-1\\\imath \\1\end{pmatrix},\quad y_2=\begin{pmatrix}-1\\1\\-1\\1\end{pmatrix},\quad y_3=\begin{pmatrix}\imath\\-1\\-\imath \\1\end{pmatrix}\right\}\]
in terms of the basis $X_1$.  We compute $R$ in basis $Y_1$ in two steps; 
\[ z_1=x_1+x_3,\quad z_2=x_1-x_3,\quad t_1=x_2+x_4,\quad t_2=x_2-x_4\]
\[ R(z_1\tens z_2)=z_2\tens z_1, \quad R(z_1\tens t_1)=t_1\tens z_1,\quad R(z_1\tens t_2)=-t_2\tens z_1\]
\[ R(z_2\tens t_1)=-t_1\tens z_2,\quad R(z_2\tens t_2)=t_2\tens z_2,\quad R(t_1\tens t_2)=t_2\tens t_1\]
and then
\[ y_0=z_1+t_1,\quad y_1=-\imath z_2-t_2,\quad y_2=-z_1+t_1,\quad y_3=\imath z_2-t_2\]
\[ R(y_1\tens y_2)=R((\imath z_2+t_2)\tens(z_1-t_1))=\imath z_1\tens z_2+\imath t_1\tens z_2-z_1\tens t_2-t_1\tens t_2\]
\[ y_2\tens y_1=(z_1-t_1)\tens (\imath z_2+t_2)=\imath z_1\tens z_2+z_1\tens t_2-\imath t_1\tens z_2-t_1\tens t_2\]
Now, by our assumption $R(y_1\tens y_2)=\mu (y_2\tens y_1)$. Here $\mu=\mu_{12}$. Hence, comparing, we need
\[ (\mu-1)(\imath z_1\tens z_2)-(\mu+1)\imath(t_1\tens z_2)+(\mu+1)(z_1\tens t_2)-(\mu-1)(t_1\tens t_2)=0\]
which requires both $\mu+1=0$ and $\mu-1=0$, which is not possible. 


\section{Structure of $\Acal(k,X,r)$}

Let $R:V\tens V\to V\tens V$ be an invertible solution of the QYBE over a field $k$. The quantum space associated to this data is the algebra $k_R[V]$ defined as a quotient of the tensor algebra $TV$ modulo the ideal generated by elements of the form $v\tens w-R(v\tens w)$. This definition is basis-independent but in the case $V=kX$ and $R$ the linear extension of a set-theoretic solution it is clear that $k_R[V]$ computed in this basis becomes the quantum algebra $\Acal(k,X,r)$. 

On the other hand it is immediate from the form of $R$ found in the main theorem above (Theorem~\ref{maintheorem}) that for a finite square-free multipermutation solution $r$ of level 2 the vector space $V=\C X$ has a new basis of generators $\{y^i_\eta\ |\ 1\leq i\leq t,\ \eta\in\Lambda_i\}$ in which the relations of $\C_R[V]$ take the diagonal form
\begin{equation}
\label{newrelations}
 y_{\eta}^i . y_{\zeta}^j=\frac{\mu^i_{j, \zeta}}{\mu^j_{i, \eta}}
y_{\zeta}^j . y_{\eta}^i,\; 1\leq i \le j \leq t, \; \eta \in  \Lambda_i, \; \zeta \in \Lambda_j,
\end{equation}
where the complex numbers $\mu_{j, \eta}^{i}$ etc., are roots of unity defined explicitly
in (\ref{theq1}). When $i=j$ we require $\eta<\zeta$ but note that due to the properties of the coefficients $\mu$ we have in in fact that
\begin{equation}\label{newrelations2}
y_{\eta}^i y_{\zeta}^i=y_{\zeta}^iy_{\eta}^i, \quad  \eta< \zeta \in \Lambda_i,\quad 1\leq i\leq t.
\end{equation}
Since the change from the $X$ to the new $Y$ basis does not change the algebra, clearly this algebra is also isomorphic to the quantum algebra  $\Acal(\mathbb C, X, r)$. This is our first immediate corollary.

In this section we give two points of view on the further particular structure of $\Acal(k,X,r)$ in the form $k_R[V]$ for our solutions of mulitpermutation level 2.   The first result concerns its structure as a cotwist indicated by the particular `ratio' form of $R$ in Theorem~\ref{maintheorem}.

\subsection{Generalities on diagonal cotwist quantum algebras}

We will use Hopf algebra methods but the Hopf algebras here will be group algebras and therefore we will give the general theory of cotwists at this level. Let $H$ be a group and $\Fcal:H\times H\to k^\times$ a 2-cocycle in the sense 
\[ \Fcal(a,bc)\Fcal(b,c)=\Fcal(a,b)\Fcal(ab,c),\quad \Fcal(e,a)=1,\quad \forall a,b,c\in H\]
(here $e$ is the group identity and one can deduce that $F(a,e)=1$ holds as well). Let $\Acal$ be a $H$-graded algebra (so $\Acal=\oplus_{a\in H}\Acal_a$ into graded components, $1\in \Acal_e$ and $\Acal_a\Acal_b\subseteq \Acal_{ab}$). Then the vector space of $A$ acquires a new associative product $\bullet$, say, with
\[ f\bullet g=\Fcal(|f|,|g|)fg\]
on homogeneous elements $f,g\in\Acal$ of degree $|f|,|g|$. The right hand side uses the old product. We obtain new algebra $\Acal^\Fcal$ covariant under the Drinfeld-cotwisted quantum group $(kH)_\Fcal$. As a Hopf algebra the latter is the same as the group algebra $kH$ but with a cotriangular structure  $\Fcal_{21}\Fcal^{-1}$ (the product here in a certain convolution algebra). This means that the natural category in which $\Acal^\Fcal$ lives is that of $H$-graded spaces but with a non-trivial symmetry (i.e. involutive brading) $\Psi$ given by a ratio of $\Fcal$ and its inverse transpose. We refer to \cite{Ma:book} for details.

Now suppose that $V$ is a vector space of dimension $n$ with basis $Y=\{y_i\}$ and $R$ a solution with diagonal entries $\rho_{ij}$ in the sense
\[ R(y_i\tens y_j)=\rho_{ij}y_j\tens y_i.\]
The indices $i,j$ in this subsection will become multi-indices $(i,\eta),(j,\zeta)$ etc in our application. We suppose further that $\rho_{ij}=F_{ij}/F_{ji}$ for some $F_{ij}$. Note that since $R$ is involutive we have $\rho_{ji}=\rho_{ij}^{-1}$  and the above form provides a natural special case of this. Let $k[Y]$ (the algebra of polynomials $k[y_1,\cdots,y_n]$ in the basis)  be graded by $H=\Z^n$, the free abelian group on  generators $g_i$, with the $H$-degree of $y_i$ being $g_i$.  

Next,  we define
\[ \Fcal(g_{i_1}\cdots g_{i_k},g_{j_1}\cdots g_{j_l})=\prod_{a=1}^k\prod_{b=1}^l F_{i_a j_b}\]
This extends to negative powers by the use of $F^{-1}_{i_aj_b}$ if one or other (but not both) of the relevant $g_{i_a}$ or $g_{j_b}$ occurs inverted in the sequence. It is easy to see that this defines a
 bicharacter on $\Z^n$ and hence in particular a 2-cocycle on it. A bicharacter on a group means $\Fcal(ab,c)=\Fcal(a,c)\Fcal(b,c)$ and $\Fcal(a,bc)=\Fcal(a,b)\Fcal(a,c)$ for any $a,b,c$ in the group.

Applying the cotwist construction now gives the algebra $k[Y]^{\Fcal}$ with the same vector space as $k[Y]$ and the new product
\begin{equation}\label{bullet}(y_{i_1}y_{i_2}\cdots y_{i_k})\bullet (y_{j_1}y_{j_2}\cdots y_{j_l})=\left(\prod_{a=1}^k\prod_{b=1}^l F_{i_a j_b}\right)  y_{i_1}y_{i_2}\cdots y_{i_k} y_{j_1}y_{j_2}\cdots y_{j_l}.\end{equation}

\begin{proposition}\label{Fquant} Suppose that $R:V\tens V\to V\tens V$ has a diagonal `ratio' form in a basis $\{y_i\}$. Then $k_R[V]$ is isomorphic to $k[Y]^{\Fcal}$, i.e. a cocycle cotwist of the commutative polynomial algebra $k[Y]$.
\end{proposition}
\proof The $n(n-1)/2$ relations of $k_R[V]$ in basis $\{y_i\}$ are of the form $y_i.y_j=\rho_{ij}y_j.y_i$ and since  the coefficients in our case are clearly invertible they allow all words on the generators to be put in standard form in terms of basis $\{y_1^{s_1}\cdots y_n^{s_n}\}$. We identify the vector space of $k_R[V]$ with that of $k[Y]^{\Fcal}$ in this way. One may then verify that the products coincide. Thus $y_i\bullet y_j=F_{ij}y_i y_j=F_{ij}y_j y_i={F_{ij}\over F_{ji}}y_j\bullet y_i$ as required in degree 2. Hence there is an algebra epimorphism $k_R[V]\to k[Y]^\Fcal$. Since their bases coincide in each degree, this is an isomorphism.  \endproof

By definition, an exterior differential calculus $(\Omega(\Acal),\extd)$ on an algebra $\Acal$ means an associative $\N$-graded algebra $\Omega$ generated  by degrees $\Omega^0=\Acal$ and $\Omega^1$, and $\extd:\Omega^m\to \Omega^{m+1}$ obeying $\extd^2=0$ and a graded Leibniz rule for its action on products. 

\begin{corollary}\label{Fcalc} $k_R[V]$ has a natural  exterior differential calculus $(\Omega_R(V),\extd)$ with basis $\extd y_i$ over $k_R[V]$ and relations
\[ \extd y_i\cdot y_j=\rho_{ij} y_j.\extd y_i,\quad\forall i,j,\quad \extd y_i \cdot \extd y_j+\rho_{ij}\extd y_j\cdot\extd y_i=0,\quad\forall i\leq j.\]
Moreover, the partial derivative operators $\del^i:k_R[V]\to k_R[V]$ defined by $\extd f=\sum_i (\del^i f)\extd y_i$ are right-handed braided derivations in the sense:
\[ \del^i(f\cdot(y_{j_1}\cdots y_{j_l}))=f\cdot \del^i(y_{j_1}\cdots y_{j_l})+  \left(\prod_{b=1}^l \rho_{i j_b}\right)(\del^i f)\cdot (y_{j_1}\cdots y_{j_l})   \]
on $f$ and homogeneous $y_{j_1}\cdots y_{j_l}$ in $k_R[V]$. The product here is the noncommutative one of $k_R[V]$.
\end{corollary}
\proof We apply the cotwisting theory above to the usual exterior algebra $\Omega(k[Y])$ generated by the commutative $y_i$ and the graded-commutative $\extd y_i$. This retains  the $\Z^n$-grading where $\extd$ is given degree 0 (so $\extd y_i$ has the same degree as $y_i$) and is therefore  an $H$-comodule algebra as before.   Cotwisting now gives $y_i\bullet \extd y_j=F_{ij}y_i\extd y_j={F_{ij}\over F_{ji}} \extd y_j\bullet y_i$ as before. Similarly for the relations between the $\extd y_i$. We define this algebra $\Omega(k[Y])^{\Fcal}$ to be the exterior calculus on  $k_R[V]=k[Y]^{\Fcal}$. The map $\extd$ is unchanged at the level of the unchanged underlying vector spaces. Finally, we identify this exterior calculus with product $\bullet$ with the abstract algebra $\Omega_R(V)$ generated by $y_i,\extd y_i$ as stated, by identifying the natural bases as before. The product in the abstract algebra is now the noncommutative one (without explicitly writing $\bullet$). The relations imply that in $\Omega_R(V)$,
\[  \extd y_i.(y_{j_1}\cdots y_{j_l})=\left(\prod_{a=1}^l\rho_{i j_a}\right) y_{j_1}\cdots y_{j_l}\extd y_i\]
from which the behaviour of the partial derivatives follows. \endproof

Another feature of $k_R[V]$ is its braided Hopf algebra structure as follows. We have explained that as a cotwist it `lives' in the category of $H$-graded spaces but with a nontrivial braiding (in our case a symmetry) built from $\rho$. The braiding on monomials is
\[ \Psi(y_{i_1}\cdots y_{i_k}\tens y_{j_1}\cdots y_{j_l})=\left(\prod_{a=1}^{k}\prod_{b=1}^l\rho_{i_a j_b}\right) y_{j_1}\cdots y_{j_l}\tens y_{i_1}\cdots y_{i_k}\]
Working in this category,  $k_R[V]$ necessarily has an additive coproduct $\underline{\Delta}y_i=y_i\tens 1+1\tens y_i$ which expresses addition on the underlying braided space just as the usual $k[Y]$ has an ordinary Hopf algebra structure with additive coproduct.  We refer to \cite[Ch 10]{Ma:book} for the general theory. (In our case since $R$ is involutive there is no other $R'$-operator as in the general theory of braided linear spaces.)

\begin{proposition} \label{Fhopf} In the context above with $k_R[V]$ viewed as a  `braided linear space' the braided derivatives $\del^i$ canonically defined by
\[ \underline\Delta f(y)=f\tens 1+ \sum_i\del^i(f)\tens y_i+\cdots\]
(where we drop terms with higher products of $y_i$ on the right) coincide with the partial derivatives above obtained by cotwisting.
\end{proposition}
\proof As a Hopf algebra in a (symmetric) braided category we extend $\underline\Delta$ to products by
\[ \underline\Delta(fg)=(\underline\Delta f)(\underline\Delta g)=\sum f\o\Psi(f\t\tens g\o)g\t\]
using the notation $\underline\Delta f=\sum f\o\tens f\t$. Expanding as stated, this becomes
\[ (f\tens 1+\del^i f\tens y_i+\cdots)(g\tens 1+\del^i g\tens y_i+\cdots)=fg\tens 1+f\del^ig\tens y_i+\del^if\Psi(y_i\tens g)+\cdots\]
which gives the same results for $\del^i(fg)$ as above. Since the values of both on generators are $\del^iy_j=\delta^i{}_j$ (the Kronecker delta-function), the partial derivatives coincide. There are explicit formulae for the braided derivatives in terms of  `braided integer matrices'  generalizing the notion of $q$-integers, see  
\cite{Ma:book}. \endproof

Similar calculus and braided Hopf algebra structures on $k_R[V]$ were found in greater generality in the 1990s in a `braided approach to noncommutative geometry' \cite{Ma:book}, without assuming that $R$ is  involutive or of cotwist form. However, in the special case of a diagonal cotwist solutiion as above the formulae are rather simpler and, moreover, we are able to obtain the constructions quite explicitly on  the underlying vector spaces rather than in terms of  arguments by generators and relations as in the literature.

All of the above apply in view of Theorem~\ref{maintheorem} to $\C_R[V]$ with $R$ obtained from a finite square-free set-theoretic solution of multipermutation level 2.  We set 
\begin{equation}\label{Fmu} F_{(i,\eta),(j,\zeta)}=\mu^i_{j,\zeta}\end{equation}
where our basis of $V$ is now written as  $\{y^i_\eta\}$, i.e. labelled by the multi-index $(i,\eta)$ in place of $i$ above. We have explained that $\C_R[V]$ is $\Acal(\C,X,r)$ written as a quadratic algebra and Proposition~\ref{Fquant} asserts that this is a $\bullet$-product `quantisation' of $\C[V]$ (the usual polynomial algebra with our basis as generators).   The relations of the differential calculus  on $\Acal(\C,X,r)$ from Corollary~\ref{Fcalc} are
\begin{equation}\label{mucalc} \extd y_{\eta}^i . y_{\zeta}^j=\frac{\mu^i_{j, \zeta}}{\mu^j_{i, \eta}}
y_{\zeta}^j. \extd y_{\eta}^i,\quad (\extd y_{\eta}^i) .(\extd  y_{\zeta}^j)+ \frac{\mu^i_{j, \eta}}{\mu^j_{i, \eta}}
(\extd y_{\zeta}^j) . (\extd y_{\eta}^i)=0,\end{equation}
(we take $i\le j$  in the second set of relations and where $i=j$ then $\eta<\zeta$). The algebra is moreover a Hopf algebra in a certain symmetric monoidal category with additive coproduct in the generators as in Proposition~\ref{Fhopf}. Note that  $\Acal(\C,X,r)$, as a semigroup algebra, also has an ordinary bialgebra structure with coproduct $\Delta x=x\tens x$ on $x\in X$. This is a different coproduct from the additive braided coproduct $\underline\Delta$ on the same algebra.

\subsection{Canonical cotwist structure of $\Acal(\C,X,r)$}

While the above generalities apply in our case, they do not make use of the specific structure of the diagonal basis $\{y^i_\eta\}$. Here we provide the same $\bullet$ product for our particular quantum algebra $\C_R[V]$ (and hence $\Acal(\C,X,r)$) as a different cotwist,  this time with  finite grading group and a canonical cocycle on it.

We let 
\[ H=\prod_{1\leq i\leq t} \Gcal_i^*\times \prod_{1\leq i\leq t} \Gcal_i=\Hcal^*\times\Hcal\]
and we let $V$ be $H$-graded as follows. We let the component $V_i$ have $\Hcal$-grade $\prod_{1\leq j\leq t}\sigma^i_j$ and define the $\Hcal^*$-grading equivalently as an action of $\Hcal$. This action is just the action of each $\Gcal_i$ on each $V_i$ (with trivial action on $V_j$ if $j\ne i$). Because the group $H$ is abelian the grading extends to $\C[V]$ (by which we mean the symmetric algebra on $V$, i.e. the free commutative algebra of polynomials in a fixed basis).

Note that these definitions are basis-independent. However, the bases where the actions of the $\Gcal_i$ are diagonal are just those where the basis elements have homogeneous grade under $\Hcal^*$ and hence under the action of $H$. This is precisely our $\{y^i_\eta\}$ basis constructed in Section~4. The basis elements here have degree
\[ |y^i_\eta|=(\chi_\eta,\gamma^i)\in H;\quad \gamma^i=\prod_{1\leq j\leq t}\sigma^i_j.\] 

Next, on each subgroup $\Gcal_i^*\times \Gcal_i$ there is a canonical 2-cocycle
\[  \Fcal( (\chi, g),(\chi',h))=\chi'(g),\quad \forall \chi,\chi'\in \Gcal_i^*,\ g,h\in\Gcal_i\]
which similarly extends for each $i$ to define a 2-cocycle $\Fcal$ on $H$. 

\begin{proposition}\label{Acotwist} Let $(X,r)$ be a finite square-free multipermutation solution of level $2$ and let $V=\C X$ be $H$-graded as above. Then the canonical 2-cocycle $\Fcal$ on $H$ induces a cotwist quantization $\C[V]^\Fcal$ isomorphic to $\C_R[V]$ and hence to $\Acal(\C,X,r)$.
\end{proposition}
\proof We compute in the homogeneous basis
\begin{eqnarray*} (y^{i_1}_{\eta_1}\cdots y^{i_k}_{\eta_k})\bullet (y^{j_1}_{\zeta_1}\cdots y^{j_l}_{\zeta_l})&=&\Fcal(\chi_{\eta_1}\gamma^{i_1}\cdots \chi_{\eta_k}\gamma^{i_k}, \chi_{\zeta_1}\gamma^{j_1}\cdots \chi_{\zeta_l}\gamma^{j_l})
 y^{i_1}_{\eta_1}\cdots y^{i_k}_{\eta_k} y^{j_1}_{\zeta_1}\cdots y^{j_l}_{\zeta_l}\\
 &=&(\chi_{\zeta_1}\cdots\chi_{\zeta_l})(\gamma^{i_1}\cdots\gamma^{i_k}) y^{i_1}_{\eta_1}\cdots y^{i_k}_{\eta_k} y^{j_1}_{\zeta_1}\cdots y^{j_l}_{\zeta_l}\\
&=&\chi_{\zeta_1}(\sigma^{i_1}_{j_1}\cdots\sigma^{i_k}_{j_1})\cdots\chi_{\zeta_l}(\sigma^{i_1}_{j_l}\cdots\sigma^{i_k}_{j_l}) y^{i_1}_{\eta_1}\cdots y^{i_k}_{\eta_k} y^{j_1}_{\zeta_1}\cdots y^{j_l}_{\zeta_l}\\
&=&\left(\prod_{a=1}^k\prod_{b=1}^l\chi_{\zeta_b}(\sigma^{i_a}_{j_b})\right)y^{i_1}_{\eta_1}\cdots y^{i_k}_{\eta_k} y^{j_1}_{\zeta_1}\cdots y^{j_l}_{\zeta_l}\end{eqnarray*}
We write expressions in the group $H$ where the various factors $\Gcal_i$ and $\Gcal_j^*$ all mutually commute. We view the product of characters in the second line here as a single character in $\Hcal^*$ but the evaluation of a factor $\chi_{\zeta_b}$ is 1 except on elements in the corresponding $\Gcal_{j_b}$. Hence writing out the products of the $\gamma^{i_a}=\sigma^{i_a}_1\cdots\sigma^{i_a}_t$ and collecting the factors living in each $\Gcal_{j_b}$ gives the third line. Using that $\chi_{\zeta_b}$ is a character gives the last expression. We recognise $\mu^{i_a}_{j_b\zeta_b}$ in the product, i.e.  the same factors $F_{(i_a,\eta_a),(j_b,\zeta_b)}$ as in the previous section and the same form of product (\ref{bullet}) in terms of such factors. Hence we again obtain the algebra $\Acal(\C,X,r)$ in the $\{y^i_\eta\}$ generators.
\endproof

These results express $\Acal(\C,X,r)$ as a kind of `quantisation' of the commutative polynomial algebra $\C[V]$ (indeed, the form of $H$ and the cocycle on it are a finite version of the similar ideas that apply to the Moyal product in quantum mechancs). We similarly find the same formulae (\ref{mucalc}) for the quantum  differential calculus and for the braiding. In the present setting, however, we learn that these constructions are all $H$-graded since we have constructed them in a symmetric monoidal category of $H$-graded spaces. 

Moreover, they are canonically defined independent of any basis. The only subtlety is that other elements such as those of our original basis $X$ will not be homogeneous. However,  we can think of the $H$-grading equivalently as a {\em coaction} of the group Hopf algebra $\C H$ where everything is extended linearly as a map   $\Delta_R:V\to V\tens \C H$ (the axioms for this map are dual to the axioms of an action $H^*\times V\to V$ and are equivalent  to such an action as $H$ here is finite). On homogeneous elements  we have $\Delta_Rv=v\tens |v|$ and we extend this linearly. We also understand the cocycle now as a map $\Fcal:\C H\tens\C H \to \C$ which is the same map as before  but extended linearly. The bullet product is defined now as
\begin{equation}\label{bulletco} f\bullet g=f\bo g\bt \Fcal(f\bt\tens g\bt);\quad \Delta_R f=\sum f\bo\tens f\bt\end{equation}
where we use an `explicit' notation for  $\Delta_Rf\in V\tens H$. Provided we linearise everything in this way, we are free to change bases.  In particular, we could now use the  original basis $X$. Then we recover the following result:
\begin{corollary}\label{LquanC} $\Acal(\C,X,r)$ can be built on the same vector space as $\C[V]$ where $V=\C X$ (i.e. the polynomial algebra on $|X|$ variables labelled by $x\in X$) with bullet product
\[  x\bullet y=x \Lcal_x(y),\quad \forall x,y\in X\]
\end{corollary}
\proof  We compute $\C[V]^\Fcal$ in the basis $X$ with the explanations above that our previous maps are extended linearly. The coaction of $\C H$ on $V=\C X$ corresponding to our previous $H$-grading has the form 
\begin{equation}\label{xcoaction}  x\mapsto {1\over N_i}\sum_{\pi\in\Gcal_i} \pi(x)\tens \left(\sum_{\eta\in\Lambda_i}\chi_\eta(\pi^{-1})\chi_\eta\right)\gamma^i.\end{equation}
Here  the $\gamma^i$  part of the grading is the same on the entire space $V_i$ hence unaffected by the change of basis. The Fourier transform relating the $\{\pi(x_{i1})\}$ basis to the $\{y^i_\eta\}$ has an inverse and we use this to compute the coaction $y^i_\eta\mapsto y^i_\eta\tens|y^i_\eta|$ in the former basis as stated in (\ref{xcoaction}). One may verify this directly from $y^i_\eta=\sum_\pi \chi_\eta(\pi^{-1})\pi(x_{i1})$ and the identity $\sum_\pi\chi_\eta(\pi)\chi(\pi^{-1})=\delta_{\chi,\chi_\eta}|\Gcal_i|$.  Armed with this coaction, for $x\in X_i$, $y\in X_j$ we have,
\[ x\bullet y={1\over N_i N_j}\sum_{\pi\in\Gcal_i\\ \rho\in\Gcal_j}\pi(x)\rho(y)\sum_{\eta\in\Lambda_i\\
\zeta\in\Lambda_j}\chi_\eta(\pi^{-1})\chi_\zeta(\rho^{-1})\Fcal(\chi_\eta\gamma^i,\chi_\zeta\gamma^j)=x\sigma^i{}_j(y)\]
after a short computation. Here $\Fcal(\chi_\eta\gamma^i,\chi_\zeta\gamma^j)=\chi_\zeta(\gamma^i)=\chi_\zeta(\sigma^i_j)$. \endproof

We remark that while it is easy to see directly that the stated formula for $\bullet$ gives the right relations $x\bullet y=({}^xy)\bullet (x^y)$ on the generators, it is not {\em apriori} obvious that $\bullet$ extends to a well-defined associative product in the same vector space as $\C[V]$ and this is what the twisting theory ensures (it follows from $\Fcal$ a 2-cocycle).  In the $X$ basis we also have
differential calculus readily computed in the same way as
\begin{equation}\label{Acalc} x\bullet \extd y=(\extd {}^xy)\bullet x^y,\quad \extd x\bullet \extd y=(\extd {}^xy)\bullet (\extd x^y)\end{equation}
with $\extd^2=0$ and $\extd$ respecting the $H$-grading, part which includes that $\extd $ commutes with the action of each restriction $\Gcal_i$ and hence with the left action. Dropping the $\bullet$ symbol, this defines the exterior differential calculus $(\Omega(\Acal(\C,X,r)),\extd)$ for our quantum space. The symmetry $\Psi$ on the generators is just $r$ itself. 

Finally, we explain how the above results extend to any field $k$. The trick is that by the time we allow the linear span of characters in $\C\Hcal^*\subseteq \C H=\C \Hcal^*\tens \C \Hcal$ we are in fact working with {\em all} functions on $\Hcal$ -- we do not need to use the basis provided by the characters and can hence dispense with the characters all together. Thus, $\C \Hcal^*=\C(\Hcal)$ in a canonical way as Hopf algebras. The latter denotes the Hopf algebra of functions on $\Hcal$ (this works for any finite group). In this second form, we can now simply work with $k(\Hcal)$ for a general field. It is the Hopf algebra dual of $k\Hcal$; we refer to \cite{Ma:book} for an introduction.  One can in fact introduce characters over $k$ but we do not need to if we use these Hopf algebra methods. Thus, instead of the group $H$ above we now work with the Hopf algebra 
\[ H=k(\Hcal)\tens k\Hcal.\]
The Hopf algebra $k(\Hcal)$ has a basis $\{\delta_h\ |\ h\in\Hcal\}$ where $\delta_h$ is the Kronecker $\delta$-function  $\delta_h(h')=1$ if $h=h'$ and otherwise 0. We explained in general above that the $\Hcal^*$-part of the grading was dual to the action of $\Hcal$ on $V=kX$ which on each $V_i=k X_i$ is given by the action of $\Gcal_i$. The $\Hcal$ part of the grading is by the element $\gamma^i$ as before. Thus the coaction $V\to V\tens H$  given simply by
\[ \Delta_R v= \sum_{\pi\in \Gcal_i}\pi.v\tens (\delta_\pi\tens \gamma^i),\quad \forall\ v\in V_i.\]
where $\delta_\pi\in k(\Gcal_i)$ is viewed in $k(\Hcal)=\tens_{1\le i\le t}k(\Gcal_i)$. This agrees with our previous formula (\ref{xcoaction}) in the case where $k=\C$ since 
\[ {1\over |\Gcal_i|}\sum_{\eta\in\Lambda_i} \chi_\eta(\pi^{-1})\chi_\eta=\delta_\pi\]
when the linear combination of characters on the left is viewed as a function on $\Gcal_i$.  

Next, as Hopf algebra 2-cocycle $\Fcal:H\tens H\to k$ we take the linear map
\begin{equation}\label{FcanH} \Fcal((f\tens h)\tens (f'\tens h'))=\epsilon(f)\epsilon(h')f'(h),\quad f,f'\in k(\Hcal),\quad h,h'\in k\Hcal\end{equation}
where the map $\epsilon(f)=f(e)$ (evaluation as the group identity) and $\epsilon(h)=1$ for all $h\in \Hcal\subset k\Hcal$. The $k$-valued linear map $\epsilon$ is called the counit of the relevant Hopf algebra. We refer to \cite[Chapter 2.3]{Ma:book} for the formal definitions and axioms obeyed by $\Fcal,\epsilon$ as linear maps. The cotwisting theory applies to produce an associative algebra with $\bullet$-product (\ref{bulletco}) but now on the vector space of $k[V]$ over a general field. 

\begin{proposition} \label{cotwistk} Let $(X,r)$ be a finite square-free multipermutation solution of level 2, and let $V=kX$ be an $H$-comodule by $\Delta_R$ as above, where $H=k(\Hcal)\tens k\Hcal$. Then the Hopf algebra 2-cocycle $\Fcal$ in (\ref{FcanH}) induces a cotwist $k[V]^\Fcal$ isomorphic to $k_R[V]$ and hence to $\Acal(k,X,r)$.
\end{proposition}
\proof The computation of the $\bullet$ product is the same as before but proceeds more cleanly without the use of characters. We compute in the basis $X$. For $x\in X_i$, $y\in X_j$, we have
\[ x\bullet y=\sum_{\pi\in \Gcal_i\atop \rho\in\Gcal_j}\pi(x)\rho(y)\Fcal((\delta_\pi\tens \gamma^i)\tens(\delta_\rho\tens\gamma^j))=\sum_{\pi\in \Gcal_i\atop \rho\in\Gcal_j}\pi(x)\rho(y)\epsilon(\delta_\pi)\epsilon(\gamma^j)\delta_\rho( \gamma^i) =x \sigma^i_j(y)\]
as before. Remember that $\delta_\rho$ is viewed as a function on $\Hcal$ but is sensitive only to the part of $\gamma^i$ that lives in $\Gcal_j$, where it has value $\delta_{\rho,\sigma^i_j}$. We then find as observed before that $x\bullet y=({}^x y)\bullet (x^y)$ as required in degree 2. Also, by construction, the algebra $k[V]^{\Fcal}$  has  the same Hilbert series as $k[V]$. As $\Acal(k,X,r)$  also has the same Hilbert series when $(X,r)$ is a symmetric set, see \cite[Thm. 1.3]{TM} and \cite[Thm. 9.7 ]{T04}, we conclude that the epimorphism  $\Acal(k,X,r)\to k[V]^{\Fcal}$ is an isomorphism. \endproof

In the same way, it is clear that we have $(\Omega(\Acal(k,X,r)),\extd)$ with the same formulae as before in the $X$ basis. One also has a monoidal symmetry $\Psi$ defined over $k$ (it is built from $\Fcal$ and its `convolution inverse', see \cite[Chapter 3]{Ma:book}) with respect to which $\Acal(k,X,r)$ is a Hopf algebra in a symmetric monoidal category.

\subsection{Structure of $\Acal(k,X,r)$ as a braided-opposite algebra}

Here we present a different point of view on the structure of $\Acal(k,X,r)$. We use Hopf algebra methods which do not require characters and in this context we prefer to work with $\Gcal(X,r)\subseteq\Hcal$ rather than the bigger group. We also then do not need to assume that $X$ is finite. Initially, let $\Gcal$ be any group.

\begin{definition} A set $X$ is a crossed $\Gcal$-set if there is a map $|\ |:X\to \Gcal$ (the grading)
and an action of $\Gcal$ on $X$ such that $|g.x|=g|x|g^{-1}$ for all $g\in \Gcal$ and $x\in X$.
\end{definition}
The notion goes back to Whitehead but is also part of the modern theory of quantum groups. The linear span $V=k X$ becomnes a `crossed $\Gcal$-module' and this is the same thing in the case of $\Gcal$  finite as $V$ a module under the quantum double $D(\Gcal)$ of the group algebra of $\Gcal$, see \cite{Ma:book} for an introduction. Its category of modules is strictly braided and so is the category of crossed $\Gcal$-modules (for any group). In algebraic terms we require $V$ to be both $\Gcal$-graded as a vector space and to enjoy an action of $\Gcal$, and the two to be compatible. The grading can be viewed  as a coaction $V\to V\tens k\Gcal$ of the Hopf algebra $k\Gcal$ which we write explicitly as $v\mapsto \Delta_R(v)=\sum v\bo\tens v\bt\in V\tens k\Gcal$, as we have seen at the end of the previous section. Let us write the action of $\Gcal$ and hence of $k\Gcal$ as usual from the left to conform with notations in the paper. Then the braiding in the category of crossed modules $\Psi^D_{V,W}:V\tens W\to W\tens V$ can be written as
\begin{equation}\label{PsiD}\Psi^D(v\tens w)=\sum v\bt. w\tens v\bo\end{equation}

\begin{lemma}\label{XGset}Let $(X,r)$ be a square-free  involutive set-theoretic  solution of  multipermutation level 2. Then $X$ is a crossed $\Gcal$-set, where $\Gcal$ is the group defined in Section~2, and  $|x|=\Lcal_x$. 
\end{lemma}
\proof If $X$ is finite, we know from Theorem~\ref{significantth} that $\Gcal$ is abelian, but the same is true even if $X$ is infinite (we will prove this in passing in  the next section, see (\ref{Lcomm})). Since the group $\Gcal$ is Abelian we require for a crossed module that the grade of an element is not changed by the action of $\Gcal$. In the finite case the grade $\Lcal_x$ is constant on each orbit $X_i$ and by definition the action of $\Gcal$ does not take us out of an orbit. One can show by direct computation that  $\Lcal_x$ does not change if one acts on $x$ by an element of $\Gcal$ even if $X$ is infinite (see (\ref{mplf})). \endproof

Hence $V$ is a crossed $\Gcal$-module and one may verify that so is $k[V]$, the polynomial algebra on any fixed basis of $V$ (it is defined invariantly as the symmetric algebra on $V$ and in that formulation does not require $V$ to be finite-dimensional).  The grade of any monomial is the product of the grades of each generator in the monomial and in this way $k[V]$ becomes $\Gcal$-graded (we need that $\Gcal$ is Abelian). The algebra also has an action of $\Gcal$ inherited from the action on $V$. In this way $k[V]$ is an algebra in the braided category of crossed $\Gcal$-modules.

Finally, for any algebra $\Acal$ over $k$ in a $k$-linear braided category there is an opposite algebra $\Acal^{\underline{op}}$ 
with `opposite product' which we denote $\bullet$ built on the same object but with
\[ \bullet= \cdot\Psi\]
where $\cdot$ denotes the original product of $\Acal$. 

\begin{proposition}\label{Abraop} Let $V=k X$ for $(X,r)$ as in the lemma above and finite, and view $k[V]$ in the braided category of crossed $\Gcal$-modules.  Then $\Acal(k,X,r)\cong k[V]^{\underline{op}}$.
\end{proposition}
\proof The opposite product in the category of crossed modules has the form
\[ f\bullet g= \cdot(f\bt.g\tens f\bo)=(f\bt.g) f\bo=f\bo (f\bt.g),\quad\forall f,g\in k[V],\]
where we use that the initial algebra is commutative and $\Psi^D$ from (\ref{PsiD}).  If we compute this
in our $X$ basis we find 
\[ x\bullet y= x\Lcal_x(y)= x\, {}^xy.\]
It is easy to see (see remarks above) that this implies $x\bullet y={}^xy\bullet x^y$ as required. The above construction ensures that the algebra $k[V]^{\underline{op}}$ is associative and is built on the same vector space as $k[V]$ in each degree as taking the opposite product does not change the  degree. At least in the case of finite $X$ it follows that this algebra is isomorphic to  $\Acal(k,X,r)$ by the same arguments as in the proof of Proposition~\ref{cotwistk}. \endproof

We do not attempt the proof with infinite $(X,r)$ here but we note $k[V]^{\underline{op}}$ is defined and associative even in this case. In the finite case, we obtain $\Acal(k,X,r)$ on the vector space of $k[V]$ with the same $\bullet$ product as in Proposition~\ref{cotwistk}  but this time we understand its product as the `braided-opposite', which is a novel approach to `quantisation' via braided categories.   Moreover, this approach works even for infinite $X$ provided $\Gcal$ remains abelian. From this point of view also, the classical exterior algebra $(\Omega(k[V]),\extd)$ is $\Gcal$-graded and has an action of $\Gcal$ by declaring that it commutes with $\extd$. In this way it too lies in the braided category of crossed $\Gcal$-modules and has a braided opposite product. The axioms of a differential exterior algebra in noncommutative geometry do not make use of any braiding (they make sense in any monoidal category with suitable additive and $k$-linear structure), hence we can view the entire construction in the braided category of crossed $\Gcal$-modules. Now, if $\Acal$ in a braided category has differential exterior algebra $(\Omega,d)$ then one can show by braid-diagram methods \cite{Ma:book} that $\Acal^{\underline{op}}$ has differential exterior algebra $(\Omega^{\underline{op}},\extd)$ where we use the graded braided opposite coproduct (the term `graded' here refers to the $\Z_2$ grading in terms of odd or even differential forms and entails an additional sign according to the degree).  Hence we can `quantise' the classical exterior algebra to obtain one on $k[V]^{\underline{op}}$ by using the graded braided-opposite product which, since the initial algebra and exterior algebra are (graded) commutative, comes out as
\[ f\bullet\extd g=f\bo \extd (f\bt. g),\quad \extd f\bullet\extd g=(\extd f\bo)\extd (f\bt. g).\]
This gives the same calculus (\ref{Acalc}) when applied to $\Acal(k,X,r)$.
Once again, these constructions are basis-independent. Over $\C$, we can compute it in the $\{y^i_\eta\}$ basis if we prefer with the same result (\ref{mucalc}). Similarly, the monoidal symmetry $\Psi$ with respect to which $\Acal(\C,X,r)$ was a Hopf algebra as in Section~6.1 can also be expressed in terms of $\Psi^D$ and its inverse if one wishes to use this approach. 

\begin{remark} Motivated by the above, we note that if we are not interested in linking up with our original diagonalisation problem,  we can also redo Proposition~\ref{cotwistk} using $\Gcal$ in place of $\Hcal$ provided $\Gcal$ is finite. We define $V$ to be a $H=k(\Gcal)\tens k\Gcal$-comodule by
\[ \Delta_R(x)=\sum_{g\in\Gcal}g(x)\tens(\delta_g\tens \Lcal_x)\]
and we define the canonical 2-cocycle $\Fcal:H\tens H\to k$ by the same formula as in (\ref{FcanH}) but now with $f,f'\in k(\Gcal)$ and $h,h'\in k\Gcal$. We obtain that $\Acal(k,X,r)$ is the cotwist of $k[V]^\Fcal$ under the coaction of $k(\Gcal)\tens k\Gcal$ by a similar computation to that in Proposition~\ref{cotwistk}.

 Proposition~\ref{Abraop} when $\Gcal$ is finite is then equivalent to this version of the cotwist by the action of $\Gcal$. Thus, a coaction of $k(\Gcal)$ on a vector space $V$ is equivalent to an action of $k\Gcal$. Hence a coaction of $k(\Gcal)\tens k\Gcal$ is equivalent to an action of $k\Gcal$ and a coaction of $k\Gcal$. In our case where $\Gcal$ is abelian these make $V$ into a crossed $k\Gcal$-module which evidently is the linearization of the crossed $\Gcal$-set in Lemma~\ref{XGset}.

Going the other way, one could similarly view the $k(\Hcal)\tens k\Hcal$ coaction in Proposition~\ref{cotwistk} as making $V$ and hence $k[V]$ into a crossed $\Hcal$-module. The coaction is given by the $\Hcal$ grade of $V_i$ as $\gamma^i$ which is just the image of the $\Gcal$-grade used above under the canonical inclusion $\Gcal\subseteq \Hcal$, while the action is the action of $\Gcal_i$ on $V_i$ which pulls back to the $\Gcal$ action under this inclusion. Therefore we necessarily obtain the same $\Psi^D$ from this crossed module as the one used in Proposition~\ref{Abraop}, and hence the same braided-opposite algebra even though the category is slightly different.  
\end{remark}

 \section{Epilogue}
 
 In the course of our study of $A(k,X,r)$ we observed in Lemma~\ref{XGset} that for $(X,r)$ a finite square-free of  multipermutation level 2 the set $X$ is a crossed $\Gcal$-set. As such it acquires a braiding of which the linearisation is $\Psi^D$ in (\ref{PsiD}). The underlying braiding, which we denote $f$, is, however, set-theoretic. Hence we have the following result  at the level of $(X,r)$ analogous to our cotwisting result for $\Acal(k,X,r)$. We can state and verify it directly in maximal generality. We let $\tau$ denote the `flip' map $\tau(x,y)=(y,x)$.

 \begin{proposition}\label{rcotwist} Let $(X,r)$ be a nondegenerate quadratic set  and define an associated non-degenerate quadratic set $(X,f)$ where
 \[ f:X\times X\to X\times X,\quad f(x,y)=({}^xy,x).\]
 Then
 \begin{enumerate}
 \item $r=f\circ\tau\circ f^{-1}$ {\em iff} $r$ is involutive with {\bf lri} (or {\em iff} $r$ is involutive and cyclic).
\item Suppose  $r$ has {\bf lri}. Then any two of the following imply the third and (1): 
 \begin{enumerate} \item  ${}^{(x^y)}z={}^xz$ for all $x,y,z\in X$ \item $f$ obeys the YBE \item $r$ obeys the YBE \end{enumerate}
  \item Suppose that $r$ is square free. Then the above hold  {\em iff} $(X,r)$ is a symmetric set of multipermutation level $\le 2$.
 \end{enumerate}
 \end{proposition}
 \proof First, note that since $r$ is nondegenerate each $\Lcal_y$ is invertible hence so is $f$, with $f^{-1}(x,y)=(y,(\Lcal_y)^{-1}x)$. The right actions associated to $f$ are by the identity map which is also invertible and its left actions are the same as for $r$ and hence each are invertible. Hence $(X,f)$ is a nondegenerate quadratic set. 
 
 For part (1), writing out $r=f\circ\tau\circ f^{-1}$, we require $(\Lcal_y)^{-1}x=x^y$, i.e. $x={}^y(x^y)$, which is equivalent to {\bf lri} by \cite[Lemma~2.19]{TSh07}, and we require ${}^{(\Lcal_y)^{-1}x}y={}^xy$ which now reduces to the {\bf cl1}.  Thus, if $r$ has the form stated, clearly it is involutive and {\bf lri} holds, while conversely if $r$ is involutive and obeys {\bf lri} then by  \cite[Prop. 2.25]{TSh07} it is cyclic and hence has the form stated. Similarly if $r$ is cyclic and involutive then by \cite[Prop. 2.25]{TSh07} {\bf lri} holds. 
 
 For part (2), we first compute directly that
 \begin{equation}\label{fYBE} {\rm YBE\ for\ }f \ \Leftrightarrow {}^x({}^yz)={}^{({}^xy)}({}^xz),\quad\forall x,y,z\in X.\end{equation}
   Hence, if {\bf l1} is assumed (for example if $r$ obeys the YBE), this condition becomes  ${}^{({}^xy)}({}^{x^y}z)={}^{({}^xy)}({}^xz)$, which by  nondegeneracy  holds if and only if the condition (a) holds (we do not in fact need to assume {\bf lri} or the full YBE in this direction). 
   
If instead we assume {\bf lri} and (a). Then {\bf l1}, {\bf l2}, {\bf lr3} (the conditions for $r$ to obey the YBE) each reduce to ${}^x({}^yz)={}^y({}^xz)$ for all $x,y,z\in X$. Similarly the above condition (\ref{fYBE}) reduces to this same commutativity condition. Hence under the assumption of {\bf lri} and (a) for $r$:
\begin{equation}\label{Lcomm} \textrm{YBE for $r$} \Leftrightarrow \textrm{$\Lcal_x\Lcal_y=\Lcal_y\Lcal_x$, $\forall\ x,y\in X$}\Leftrightarrow\textrm{YBE for $f$}.\end{equation}
Moreover, (2)(a) implies that $(X,r)$ is cyclic and hence if  {\bf lri} also holds then part (1) applies.

For (3), the argument is similar to the proof in Proposition~\ref{sigprop}, but we will be careful now to allow $X$ infinite. Thus, condition $mpl(X,r)\le 2$ means $\Lcal_{[y]}=\Lcal_{[y']}$ for all $y,y'$ (so that the second retract is trivial), i.e. ${}^{[y]}[x]={}^{[y']}[x]$ for all $y,y',x$. Equality here is in $[X]$, i.e. 
${}^{({}^y x)}z={}^{({}^{y'} x)}z$ for all $y,y',x,z$. If $(X,r)$ is square-free then this is equivalent to ${}^{({}^yx)}z={}^xz$ (in one direction, set $y'=x$) as in (\ref{mp2}).  Hence under the assumption of {\bf lri} and  square-free,  
\begin{equation}\label{mplf} mpl(X,r)\le 2\quad\Leftrightarrow {}^{x^y}z={}^xz,\quad \forall x,y,z\in X.\end{equation}
In this case, if we assume (2) then by (2)(c) $(X,r)$ also obeys the YBE and by (1) it is involutive.  Conversely, if $(X,r)$ is a square-free involutive solution of the YBE then it is necessarily {\bf lri} by \cite[Cor. 2.31]{TSh07} and if   $mpl(X,r)\le 2$ then by the above, (2)(a) holds as well as (2)(c). Hence all conditions hold.  Note that in this case the group $\Gcal(X,r)$ is abelian by (\ref{Lcomm}) without assuming that $(X,r)$ is finite as in Theorem~\ref{significantth}(1).  \endproof

\subsection*{Acknowledgements} The first author would like warmly to thank 
Jos\'e G\'omez-Torrecillas   for suggesting in 2006 the problem  of investigating when $\Acal(\C,X,r)$ has
a diagonal form of relations. She would also like warmly to thank Peter Cameron
 for valuable and fruitful discussions and 
  for his beautiful  explanation of the elementary theory of characters, essential to the present
 version of the paper. The work was completed while the first author was attending the Isaac
Newton Institute Programme on Combinatorics and Statistical Mechanics 
(CSM) 2008, she thanks the Isaac Newton Institute for local support and
for the inspiring working atmosphere. Finally,  she thanks the ICTP in Trieste for support during the first stage of the project.

\end{document}